\documentclass[11pt,reqno]{amsart}
\usepackage{amssymb,amscd,amsbsy}
\usepackage{amsthm}
\newcommand{\lb}{\linebreak}

\renewcommand{\a}{\alpha}
\renewcommand{\b}{\beta}

\newcommand{\z}{\zeta}

\renewcommand{\l}{\lambda}

\newcommand{\s}{\sigma}
\renewcommand{\t}{\tau}

\newcommand{\f}{\varphi}
\renewcommand{\o}{\omega}

\newcommand{\D}{\Delta}
\renewcommand{\L}{\Lambda}

\renewcommand{\O}{\Omega}

\newcommand{\B}{{\mathcal B}}
\newcommand{\cc}{{\mathcal C}}
\newcommand{\E}{{\mathcal E}}
\newcommand{\cd}{{\mathcal D}}
\newcommand{\F}{{\mathcal F}}

\newcommand{\h}{{\mathcal H}}

\newcommand{\K}{{\mathcal K}}
\newcommand{\cL}{{\mathcal L}}

\newcommand{\X}{{\mathcal X}}
\newcommand{\Y}{{\mathcal Y}}

\newcommand{\V}{{\mathcal V}}

\newcommand{\T}{{\Bbb T}}
\newcommand{\pp}{{\Bbb P}}
\newcommand{\dd}{{\Bbb D}}
\newcommand{\R}{{\Bbb R}}
\newcommand{\Z}{{\Bbb Z}}

\newcommand{\0}{{\boldsymbol{0}}}

\newcommand{\bs}{\boldsymbol}

\newcommand{\m}{{\boldsymbol m}}
\newcommand{\bS}{{\boldsymbol S}}

\newcommand{\rf}[1]{(\ref{#1})}

\newcommand{\df}{\stackrel{\mathrm{def}}{=}}
\newcommand{\dist}{\operatorname{dist}}

\newcommand{\spn}{\operatorname{span}}
\newcommand{\supp}{\operatorname{supp}}
\newcommand{\clos}{\operatorname{clos}}

\newcommand{\trace}{\operatorname{trace}}

\newcommand{\const}{\operatorname{const}}

\newcommand{\eeq}{\end{equation}}
\newcommand{\beq}{\begin{equation}}
\newcommand{\bay}{\begin{eqnarray}}
\newcommand{\ba}{\begin{align*}}
\newcommand{\ea}{\end{align*}}
\newcommand{\ey}{\end{eqnarray}}
\newcommand{\bey}{\begin{eqnarray*}}
\newcommand{\eey}{\end{eqnarray*}}

\newcommand{\imp}{\Rightarrow}
\newcommand{\be}{\infty}

\newcommand{\bl}{\blacksquare}

\newcommand{\Pf}{{\bf Proof. }}

\newcommand{\ov}{\overline}

\newtheorem{thm}{\hspace{\parindent}Theorem}[section]

\newtheorem{cor}[thm]{\hspace{\parindent}Corollary}
\newtheorem{lem}[thm]{\hspace{\parindent}Lemma}

\theoremstyle{remark}

\newtheorem*{rem*}{Remark}

\newcommand\up{\upsilon}
\newcommand\Li{{\rm Lip}}
\newcommand\fM{\frak M}
\newcommand\cZ{\mathcal{Z}}
\newcommand\dg{\frak D}

\begin{document}
\hfill To the memory of Paul Halmos

\

\

\newcommand{\vse}{\vspace{.2in}}
\numberwithin{equation}{section}

\title{The behavior of functions of operators\\ under perturbations}
\author{V.V. Peller}
\thanks{The author is partially supported by NSF grant DMS 0700995}

\begin{abstract}
This is a survey article. We consider different problems in connection with the behavior of functions of operators under perturbations of operators. We deal with three classes of operators: unitary operators, self-adjoint operators, and contractions.
We study operator Lipschitz and operator differentiable functions. We also study the behavior of functions under perturbations of an operator by an operator of Schatten--von Neumann class $\bS_p$ and apply the results to the Livschits--Krein and Koplienko--Neidhardt
trace formulae. We also include in this survey article recent unexpected results obtained in a joint paper with Aleksandrov on
operator H\"older--Zygmund functions. 
\end{abstract}

\maketitle

\newcommand{\mt}{{\mathcal T}}

\

\begin{center}
{\Large Contents}
\end{center}

\

\begin{enumerate}
\item[1.] Introduction \quad\dotfill \pageref{s1}
\item[2.]  Double operator integrals \quad\dotfill \pageref{s2}
\item[3.] Multiple operator integrals  \quad\dotfill \pageref{s3}
\item[4.] Besov spaces \quad\dotfill \pageref{s4}
\item[5.] Nuclearity of Hankel operators \quad\dotfill \pageref{s5}
\item[6.]  Operator Lipschitz and operator differentiable functions. Sufficient \\ conditions \quad\dotfill \pageref{s6}
\item[7.]  Operator Lipschitz and operator differentiable functions. Necessary \\ conditions \quad\dotfill \pageref{s7}
\item[8.]  Higher order operator derivatives \quad\dotfill \pageref{s8}
\item[9.]  The case of contractions\quad\dotfill \pageref{s9}
\item[10.]  Operator H\"older--Zygmund functions\quad\dotfill \pageref{s10}
\item[11.]  Livshits--Krein trace formulae\quad\dotfill \pageref{s11}
\item[12.]  Koplienko--Neidhardt trace formulae\quad\dotfill \pageref{s12}
\item[13.]  Perturbations of class $\bS_p$\quad\dotfill \pageref{s13}
\item[] References \quad\dotfill \pageref{bibl}
\end{enumerate}

\

\

\setcounter{section}{0}
\section{\bf Introduction}
\setcounter{equation}{0}
\label{s1}

\

This survey article is devoted to problems in perturbation theory that arise in an attempt to understand the behavior of the function $f(A)$ of an operator $A$ under perturbations of $A$. 

Consider the following example of such problems. Suppose that $\f$ is a function on the real line $\R$, $A$ is a self-adjoint operator on Hilbert space. The spectral theorem for self-adjoint operators allows us to define the function $\f(A)$ of $A$. Suppose that $K$ is a bounded self-adjoint operator. We can ask the question of when the function
\bay
\label{dif}
K\mapsto\f(A+K)
\ey
is differentiable. We can consider differentiability in the sense of G\^ateaux or in the sense of Fr\'echet and we can consider the problem for bounded self-adjoint operators $A$ or for arbitrary self-adjoint operators (i.e., not necessarily bounded). 

It is obvious that for this map to be differentiable (in the sense of Fr\'echet) it is necessary that $\f$ is a differentiable function on $\R$. Functions, for which the map
\rf{dif} is differentiable are called {\it operator differentiable}. This term needs a clarification: we can consider operator differentiable functions in the sense of G\^ateaux or Fr\'echet and we can consider this property for bounded $A$ or arbitrary self-adjoint operators $A$.
In \cite{W} Widom asked the question: when are differentiable functions differentiable?

We also consider in this survey the problem of the existence of higher operator derivatives.

Another example of problems of perturbation theory we are going to consider in this survey is the problem to describe {\it operator Lipschitz functions}, i.e., functions $\f$ on $\R$, for which
\bay
\label{opL}
\|\f(A)-\f(B)\|\le\const\|A-B\|
\ey
for self-adjoint operators $A$ and $B$. Sometimes such functions are called uniformly operator Lipschitz.
Here $A$ and $B$ are allowed to  be unbounded provided the difference $A-B$ is bounded.
If $\f$ is a function, for which \rf{opL} holds for bounded operators $A$ and $B$ with a constant that can depend on $\|A\|$ and $\|B\|$,
then $\f$ is called {\it locally operator Lipschitz}.
It is easy to see that if $\f$ is operator Lipschitz, then $\f$ must be a Lipschitz function, i.e.,
\bay
\label{Lip}
|\f(x)-\f(y)|\le\const|x-y|,\quad x,\,y\in\R,
\ey
and if $\f$ is locally operator Lipschitz, then $\f$ is locally a Lipschitz function, i.e., \rf{Lip} must hold on each bounded subset of $\R$.

We also consider in this survey the problem for which functions $\f$
\bay
\label{oHo}
\|\f(A)-\f(B)\|\le\const\|A-B\|^\a
\ey
for self-adjoint operators $A$ and $B$. Here $0<\a<1$. 
If $\f$ satisfies \rf{oHo}, it is called
 {\it an operator H\"older function of order} $\a$. 
 Again, it is obvious that for $\f$ to be operator H\"older of order $\a$ it is necessary that $\f$ belongs to the H\"older class
$\L_\a(\R)$, i.e.,
\bay
\label{Hol}
|\f(x)-\f(y)|\le\const|x-y|^\a
\ey
and if $\f$ is operator H\"older of order $\a$, then $\f$ must satisfy \rf{Hol} on each bounded subset of $\R$.

We also consider functions $f$ on $\R$ for which
\bay
\label{oZy}
\|\f(A-K)-2\f(A)+\f(A+K)\|\le\const\|K\|
\ey
for selfadjoint operators $A$ and $K$. Functions $\f$ satisfying \rf{oZy} are called {\it operator Zygmund functions}. 

In this paper we also study the whole scale of operator H\"older--Zygmund classes.

Another group of problems we are going to consider is the behavior of functions of operators under perturbations of trace class (or other classes of operators). In particular, the problem to describe the class of functions $f$, for which 
$$
f(A+K)-f(A)\in\bS_1\quad\mbox{whenever}\quad K\in\bS_1,
$$
is very important in connection with the Livshits--Krein trace formula. We use the notation $\bS_p$ for Schatten--von Neumann classes.

We also consider problems of perturbation theory related to the Koplienko trace formula, which deals with Hilbert--Schmidt perturbations.

It is also important to study similar problems for unitary operators and functions on the unit circle $\T$ and for contractions and 
analytic functions in the unit disk $\dd$.

The study of the problem of differentiability of functions of self-adjoint operators on Hilbert space was initiated 
By Daletskii and S.G. Krein in \cite{DK}. They showed that for a function $\f$ on the real line $\R$ of class $C^2$ and for 
bounded self-adjoint operators $A$ and $K$ the function
\bay
\label{sso}
t\mapsto \f(A+tK)
\ey
is differentiable in the operator norm and the derivative can be computed in terms of double operator integrals:
\bay
\label{DKF}
\frac{d}{dt}\f(A+tK)\Big|_{t=0}=\iint\limits_{\R\times\R}\frac{\f(x)-\f(y)}{x-y}\,dE_A(x)\,K\,dE_A(y),
\ey
where $E_A$ is the spectral measure of $A$. The expression on the right is a double operator integral. The beautiful theory of double operator integrals due to Birman and Solomyak was created later in \cite{BS1}, \cite{BS2}, and \cite{BS3} (see also 
the survey article \cite{BS5}). A brief introduction into the theory of double operator integrals will be given in \S\,2.

The condition $\f\in C^2$ was relaxed by Birman and Solomyak in  \cite{BS3}: they proved that the function \rf{sso} is differentiable and the Daletskii--Krein formula \rf{DKF} holds under the condition that $\f$ is differentiable and the derivative $\f'$ satisfies a H\"older condition of order $\a$ for some $\a>0$. The approach of Birman and Solomyak is based on their formula
\bay
\label{BSdof}
\f(A+K)-\f(A)=\iint\limits_{\R\times\R}\frac{\f(x)-\f(y)}{x-y}\,dE_{A+K}(x)\,B\,dE_A(y).
\ey
Actually, Birman and Solomyak showed in \cite{BS3} that formula \rf{BSdof}
is valid under the condition that the divided difference $\dg \f$,
$$
(\dg \f)(x,y)=\frac{\f(x)-\f(y)}{x-y},
$$
 is a Schur multiplier of the space of all bounded linear operators
(see \S\,2 for definitions).

Nevertheless, Farforovskaya proved in  \cite{F0} that the condition $\f\in C^1$ does not imply that $\f$ is operator Lipschitz, which implies
that  the condition $\f\in C^1$ is not sufficient for the differentiability of the map \rf{sso} (see also \cite{F} and \cite{Ka}).

A further improvement was obtained in \cite{Pe1} and \cite{Pe4}: it was shown that the function \rf{sso} is differentiable and
\rf{DKF} holds under the assumption that $\f$ belongs to the Besov space $B_{\be1}^1(\R)$ (see \S\,4) and under the same
assumption $\f$ must be uniformly operator Lipschitz. In the same paper \cite{Pe1} a necessary condition was found: $\f$ must locally belong to the Besov space $B_1^1(\R)=B^1_{11}(\R)$. This necessary condition also implies that the condition $\f\in C^1$ is not sufficient. Actually, in \cite{Pe1} and \cite{Pe4} a stronger necessary condition was also obtained; see \S 7 for further discussions. Finally, we mention another sufficient condition obtained in \cite{ABF} which is slightly better than the condition \lb$\f\in B_{\be1}^1(\R)$, though I believe it is more convenient to work with Besov spaces. We refer the reader to Sections 6 and 7 of this survey for a detailed discussion.

After it had become clear that Lipschitz functions do not have to be operator Lipschitz, many mathematicians believed that
H\"older functions of order $\a$, $0<\a<1$, are not necessarily operator H\"older functions of order $\a$. In \cite{F0} the following upper estimate for self-adjoint operators $A$ and $B$ with spectra in an interval $[a,b]$ was obtained:
$$
\|\f(A)-\f(B)\|\le\const\|\f\|_{\L_\a(\R)}\left(\log\left(\frac{b-a}{\|A-B\|}+1\right)+1\right)^2\|A-B\|^\a,
$$
where $\f\in\L_\a(\R)$. A similar inequality was obtained in \cite{FN} for arbitrary moduli of continuity.

{\it Surprisingly, it turns out that the logarithmic factor in the above inequality is unnecessary. In other words, for an arbitrary $\a\in(0,1)$,
a H\"older function of order $\a$ must be operator H\"older of order $\a$. Moreover, the same is true for Zygmund functions and for the whole scale of H\"older--Zygmund classes.} This has been proved recently in \cite{AP1}, see also \cite{AP}. We discuss the results of \cite{AP1} in \S\,10.

The problem of the existence of higher order derivatives of the function \rf{sso} was studied in \cite{S} where it was shown that under certain assumptions on $\f$, the function \rf{sso} has a second derivative that can be expressed in terms of 
the following triple operator integral:
$$
\frac{d^2}{dt^2}\f(A+tB)\Big|_{t=0}=2\iiint\limits_{\R\times\R\times\R}\left(\dg^2\f\right)(x,y,z)\,dE_A(x)\,B\,dE_A(y)\,B\,dE_A(z),
$$
where $\dg^2\f$ stands for the divided difference of order 2 (see \S\,8 for the definition).
To interpret triple operator integrals, repeated integration was used in \cite{S}
(see also the earlier paper \cite{Pa}, in which an attempt to define multiple operator integrals was given). 
However, the class of integrable functions in \cite{Pa} and \cite{S} was rather narrow and the assumption on $\f$ imposed in \cite{S} for the existence of the second operator derivative was too restrictive.
Similar results were also obtained in \cite{S} for the $n$th derivative and multiple operator integrals.

In \cite{Pe8} a new approach to multiple operator integrals was given. It is based on integral projective tensor products of $L^\be$ spaces and gives a much broader class of integrable functions than under the approaches of \cite{Pa} and \cite{S}. It was shown in \cite{Pe8} that under the assumption that $f$ belongs to the Besov space $B_{\be1}^n(\R)$ the function \rf{sso} has $n$ derivatives and the $n$th derivative can be expressed in terms of a multiple operator integral. Similar results were also obtained in \cite{Pe8} in the case of an unbounded self-adjoint operator $A$.

To study the problem of differentiability of functions of unitary operators, we should consider a Borel function $f$ on the unit circle
$\T$ and the map
$$
U\mapsto f(U),
$$
where $U$ is a unitary operator on Hilbert space. If $U$ and $V$ are unitary operators and $V=e^{{\rm i}A}U$, 
where $A$ is a self-adjoint operator, we can consider the one-parameter family of unitary operators
$$
e^{{\rm i}tA}U,\quad0\le t\le1,
$$
and study the question of the differentiability of the function 
$$
t\mapsto \big(e^{{\rm i}tA}U\big)
$$
and the question of the existence of its higher derivatives.
The results in the case of unitary operators are similar to the results for self-adjoint operators, see \cite{BS3}, \cite{Pe1},
\cite{ABF}, \cite{Pe8}.

Similar questions can be also considered for functions of contractions. It turns out that to study such problems for contractions, 
one can use double and multiple operator integrals with respect to semi-spectral measures. This approach was proposed in \cite{Pe2}
and \cite{Pe9}. We discuss these issues in \S\,9.

In Sections 2 and 3 of this survey we give a brief introduction in double operator integrals and multiple operator integrals.
In \S 4 we introduce the reader to Besov spaces.

In \S 5 we define Hankel operators and state the nuclearity criterion obtained in \cite{Pe0}. It turns out that Hankel operators play an important role in perturbation theory and the nuclearity criterion is used in \S\,7 to obtain necessary conditions for operator differentiability and operator Lipschitzness.

The last 3 sections are devoted to perturbations of operators by operators from Schatten--von Neumann classes $\bS_p$.
In \S 11 we discuss the problem of classifying
the functions $\f$ for which the Livshits--Krein trace formulae are valid. This problem is closely related to the problem of classifying
the functions $\f$ on $\R$ for which  a trace class perturbation of a self-adjoint operator $A$ leads to a trace class perturbation of $\f(A)$.
We present in \S\,11 sufficient conditions obtained in \cite{Pe1} and \cite{Pe4} that improve earlier results by M.G. Krein and Birman--Solomyak. We also discuss necessary conditions.

Section 12 deals with perturbations of class $\bS_2$ and trace formulae by Koplienko and Neidhardt. We discuss the results of \cite{Pe7} that improve earlier results by Koplienko and Neidhardt.

Finally, in the last section we present recent results of \cite{AP2} (see also \cite{AP}) that concern the following problem. Suppose that $A$ and $B$ are self-adjoint operators such that $A-B\in\bS_p$ and let $\f\in\L_\a(\R)$. What can we say about $\f(A)-\f(B)$? We also discuss a similar problem for higher order differences.

It was certainly impossible to give proofs of all the results discussed in this survey. I tried to give proofs of certain key results that demonstrate principal ideas.

\

\section{\bf Double operator integrals}
\setcounter{equation}{0}
\label{s2}

\

In this section we give a brief introduction in the theory of double operator integrals
developed by Birman and Solomyak in \cite{BS1}, \cite{BS2}, and \cite{BS3}, see also their survey \cite{BS5}.

Let $(\X,E_1)$ and $(\Y,E_2)$ be spaces with spectral measures $E_1$ and $E_2$
on a Hilbert spaces $\h_1$ and $\h_2$. Let us first define double operator integrals
\bay
\label{doi}
\int\limits_{\X}\int\limits_{\Y}\Phi(x,y)\,d E_1(x)\,Q\,dE_2(y),
\ey
for bounded measurable functions $\Phi$ and operators $Q:\h_2\to\h_1$
of Hilbert--Schmidt class $\bS_2$. Consider the set function $F$ whose values are orthogonal
projections on the Hilbert space $\bS_2(\h_2,\h_1)$ of Hilbert--Schmidt operators from $\h_2$ to $\h_1$, which is defined 
on measurable rectangles by
$$
F(\D_1\times\D_2)Q=E_1(\D_1)QE_2(\D_2),\quad Q\in\bS_2(\h_2,\h_1),
$$ 
$\D_1$ and $\D_2$ being measurable subsets of $\X$ and $\Y$. Note that left multiplication by $E_1(\D_1)$
obviously commutes with right multiplication by $E_2(\D_2)$.

 It was shown in \cite{BS4} that $F$ extends to a spectral measure on 
$\X\times\Y$. If $\Phi$ is a bounded measurable function on $\X\times\Y$, we define
$$
\int\limits_{\X}\int\limits_{\Y}\Phi(x,y)\,d E_1(x)\,Q\,dE_2(y)=
\left(\,\,\int\limits_{\X_1\times\X_2}\Phi\,dF\right)Q.
$$
Clearly,
$$
\left\|\,\,\int\limits_{\X}\int\limits_{\Y}\Phi(x,y)\,dE_1(x)\,Q\,dE_2(y)\right\|_{\bS_2}
\le\|\Phi\|_{L^\be}\|Q\|_{\bS_2}.
$$

%It is easy to see from the definition of double operator integrals in the case $Q\in\bS_2$ 
%that if $\{\Phi_n\}_{n\ge1}$ is a sequence of measurable functions such that 
%$$
%\lim_{n\to\be}\Phi_n(\l,\mu)=\Phi(\l,\mu),\quad \l\in\X_1,~\mu\in\X_2,
%$$
%and
%$$
%\sup_n\sup_{\l,\,\mu\in\X_1\times\X_2}|\Phi_n(\l,\mu)|<\be,
%$$
%then
%\bay
%\label{poto}
%\lim_{n\to\be}\left\|~\int\limits_{\X_1}\int\limits_{\X_2}\Phi_n(\l,\mu)\,d E_1(\l)\,Q\,dE_2(\mu)
%-\int\limits_{\X_1}\int\limits_{\X_2}\Phi(\l,\mu)\,d E_1(\l)\,Q\,dE_2(\mu)\right\|_{\bS_2}=0.
%\ey

If the transformer
$$
Q\mapsto\int\limits_{\X}\int\limits_{\Y}\Phi(x,y)\,d E_1(x)\,Q\,dE_2(y)
$$
maps the trace class $\bS_1$ into itself, we say that $\Phi$ is a {\it Schur multiplier of $\bS_1$ associated with 
the spectral measures $E_1$ and $E_2$}. In
this case the transformer
\bay
\label{tra}
Q\mapsto\int\limits_{\Y}\int\limits_{\X}\Phi(x,y)\,d E_2(y)\,Q\,dE_1(x),\quad Q\in \bS_2(\h_1,\h_2),
\ey
extends by duality to a bounded linear transformer on the space of bounded linear operators from $\h_1$ to $\h_2$
and we say that the function $\Psi$ on $\X_2\times\X_1$ defined by 
$$
\Psi(y,x)=\Phi(x,y)
$$
is {\it a Schur multiplier of the space of bounded linear operators associated with $E_2$ and $E_1$}.
We denote the space of such Schur multipliers by $\fM(E_2,E_1)$

%We denote by $\fM(E_1,E_2)$ the space of Schur multipliers of $\bS_1$ associated with
%the spectral measures $E$ and $F$. The norm of $\psi$ in $\fM(E,F)$ is, by definition, the norm of the
%transformer \rf{tra} on the space of bounded linear operators.

Birman in Solomyak obtained in \cite{BS3} the following result:

\begin{thm}
\label{BST}
Let $A$ be a self-adjoint operator (not necessarily bounded) and let
$K$ be a bounded self-adjoint operator. Suppose that 
$\f$ is a continuously differentiable 
function on $\R$ such that the divided difference $\dg\f\in\fM(E_{A+K},E_A)$. Then
\bay
\label{BSF}
\f(A+K)-\f(A)=\iint\limits_{\R\times\R}\frac{\f(x)-\f(y)}{x-y}\,dE_{A+K}(x)K\,dE_A(y)
\ey
and
$$
\|\f(A+K)-\f(A)\|\le\const\|\dg \f\|_{\fM}\|K\|,
$$
where $\|\dg \f\|_{\fM}$ is the norm of $\dg \f$ in $\fM(E_{A+K},E_A)$. 
\end{thm}

In the case when $K$ belongs to the Hilbert Schmidt class $\bS_2$, the same result was established in \cite{BS3}
under weaker assumptions on $\f$:

\begin{thm}
\label{BST2}
Let $A$ be a self-adjoint operator (not necessarily bounded) and let
$K$ be a self-adjoint operator of class $\bS_2$. If $\f$ is a Lipschitz function on $\R$,
then {\em\rf{BSF}} holds,
$$
\f(A+K)-\f(A)\in\bS_2,
$$
and
$$
\|\f(A+K)-\f(A)\|_{\bS_2}\le\sup_{x\ne y}\frac{|\f(x)-\f(y)|}{|x-y|}\,\|A-B\|_{\bS_2}.
$$
\end{thm}

Note that  if $\f$ is not differentiable, $\dg \f$ is not defined on the diagonal of $\R\times\R$, but formula
\rf{BSF} still holds if we define $\dg \f$ to be zero on the diagonal.

Similar results also hold for functions on the unit circle and for unitary operators.

It is easy to see that if a function $\Phi$ on $\X\times\Y$ belongs to the {\it projective tensor
product}
$L^\be(E_1)\hat\otimes L^\be(E_2)$ of $L^\be(E_1)$ and $L^\be(E_2)$ (i.e., $\Phi$ admits a representation
\bay
\label{ptp}
\Phi(x,y)=\sum_{n\ge0}f_n(x)g_n(y),
\ey
where $f_n\in L^\be(E_1)$, $g_n\in L^\be(E_2)$, and
\bay
\label{ptpn}
\sum_{n\ge0}\|f_n\|_{L^\be}\|g_n\|_{L^\be}<\be),
\ey
then $\Phi\in\fM(E_1,E_2)$, i.e., $\Phi$ is a Schur multiplier of the space of bounded linear operators. For such functions $\Phi$ we have
$$
\int\limits_\X\int\limits_\Y\Phi(x,y)\,d E_1(x)Q\,dE_2(y)=
\sum_{n\ge0}\left(\,\int\limits_\X f_n\,dE_1\right)Q\left(\,\int\limits_\Y g_n\,dE_2\right).
$$ 
Note that if $\Phi$ belongs to the projective tensor product $L^\be(E_1)\hat\otimes L^\be(E_2)$, its norm in $L^\be(E_1)\hat\otimes L^\be(E_2)$
is, by definition, the infimum of the  left-hand side of \rf{ptpn} over all representaions \rf{ptp}.

One can define in the same way projective tensor products of other function spaces.

More generally, $\Phi$ is a Schur multiplier  if $\Phi$ 
belongs to the {\it integral projective tensor product} $L^\be(E_1)\hat\otimes_{\rm i}L^\be(E_2)$ of $L^\be(E_1)$ and $L^\be(E_2)$, i.e., $\Phi$ admits a representation
\bay
\label{ipt}
\Phi(x,y)=\int_\O f(x,\o)g(y,\o)\,d\s(\o),
\ey
where $(\O,\s)$ is a measure space, $f$ is a measurable function on $\X\times \O$,
$g$ is a measurable function on $\Y\times \O$, and
\bay
\label{ir}
\int_\O\|f(\cdot,\o)\|_{L^\be(E_1)}\|g(\cdot,\o)\|_{L^\be(E_2)}\,d\s(\o)<\be.
\ey
If $\Phi\in L^\be(E_1)\hat\otimes_{\rm i}L^\be(E_2)$, then
$$
\int\limits_\X\int\limits_\Y\Phi(x,y)\,d E_1(x)\,Q\,dE_2(y)=
\int\limits_\O\left(\,\int\limits_\X f(x,\o)\,dE_1(x)\right)Q
\left(\,\int\limits_\Y g(y,\o)\,dE_2(y)\right)\,d\s(\o).
$$
Clearly, the function 
$$
\o\mapsto \left(\,\int\limits_\X f(x,\o)\,dE_1(x)\right)Q
\left(\,\int\limits_\Y g(y,\o)\,dE_2(y)\right)
$$
is weakly measurable and
$$
\int\limits_\O\left\|\left(\,\int\limits_\X f(x,\o)\,dE_1(x)\right)T
\left(\,\int\limits_\Y g(y,\o)\,dE_2(y)\right)\right\|\,d\s(\o)<\be.
$$

Moreover, it can easily be seen that such functions $\Phi$ are Schur multipliers of an arbitrary symmetrically normed ideal of operators.

It turns out that all Schur multipliers of the space of bounded linear operators can be obtained in this way. More precisely, the following result holds (see \cite{Pe1}):

\begin{thm}
\label{SM}
 Let $\Phi$ be a measurable function on 
$\X\times\Y$. The following are equivalent:

{\rm (i)} $\Phi\in\fM(E_1,E_2)$;

{\rm (ii)} $\Phi\in L^\be(E_1)\hat\otimes_{\rm i}L^\be(E_2)$;

{\rm (iii)} there exist measurable functions $f$ on $\X\times\O$ and $g$ on $\Y\times\O$ such that
{\em\rf{ipt}} holds and
\bay
\label{bs}
\left\|\int_\O|f(\cdot,\o)|^2\,d\s(\o)\right\|_{L^\be(E_1)}
\left\|\int_\O|g(\cdot,\o)|^2\,d\s(\o)\right\|_{L^\be(E_2)}<\be.
\ey
\end{thm}

Note that the implication (iii)$\imp$(ii) was established in \cite{BS3}. Note also that
 the equivalence of (i) and (ii) is deduced from Grothendieck's theorem.
In the case of matrix Schur multipliers (this corresponds to discrete spectral measures
of multiplicity 1), the equivalence of (i) and (ii) was proved in \cite{Be}.

It is interesting to observe that if $f$ and $g$ satisfy \rf{ir}, then they also satisfy
\rf{bs}, but the converse is false. However, if $\Phi$ admits a representation of the form \rf{ipt}
with $f$ and $g$ satisfying \rf{bs}, then it also admits a (possibly different) representation of the
form \rf{ipt} with $f$ and $g$ satisfying \rf{ir}.

Let us also mention one more observation by Birman and Solomyak, see \cite{BS3}. Suppose that $\mu$ and $\nu$ are scalar measures on $\X$ and $\Y$ that are mutually absolutely continuous with 
$E_1$ and $E_2$. Let $\frak N_{\mu,\nu}$ be the class of measurable functions $k$ on $\X\times\Y$ such that the integral operator from 
$L^2(\mu)$ to $L^2(\nu)$ with kernel function $k$ belongs to the trace class $\bS_1$.

\begin{thm}
\label{zam}
A measurable function $\Phi$ on $\X\times\Y$ is a Schur multiplier of $\bS_1$ associated with $E_1$ and $E_2$
if and only if
$\Phi$ is a multiplier of $\frak N_{\mu,\nu}$, i.e.,
$$
k\in{\frak N}_{\mu,\nu}\quad\imp\quad\Phi k\in{\frak N}_{\mu,\nu}.
$$
\end{thm}

%Note that in a similar way we can define the {\it projective tensor product} $A\hat\otimes B$
%and the {\it integral projective tensor product} $A\hat\otimes_{\rm i} B$ of 
%arbitrary Banach functions spaces $A$ and $B$.

%Similar results also hold in the case of unitary operators.

\

\section{\bf Multiple operator integrals}
\setcounter{equation}{0}
\label{s3}

\

The equivalence of (i) and (ii) in the Theorem \ref{SM} suggests the idea 
explored in  \cite{Pe8} to define multiple operator integrals. 

To simplify the notation, we consider here the case of triple operator integrals; the case of arbitrary multiple operator integrals can be treated in the same way.

Let $(\X,E_1)$, $(\Y,E_2)$, and $(\cZ,E_3)$
be spaces with spectral measures $E_1$, $E_2$, and $E_3$ on Hilbert spaces $\h_1$, $\h_2$, and $\h_3$. Suppose that
$\Phi$ belongs to the integral projective tensor product
$L^\be(E_1)\hat\otimes_{\rm i}L^\be(E_2)\hat\otimes_{\rm i}L^\be(E_3)$, i.e., $\Phi$ admits a representation
\bay
\label{ttp}
\Phi(x,y,z)=\int_\O f(x,\o)g(y,\o)h(z,\o)\,d\s(\o),
\ey
where $(\O,\s)$ is a measure space, $f$ is a measurable function on $\X\times \O$,
$g$ is a measurable function on $\Y\times \O$, $h$ is a measurable function on $\cZ\times \O$,
and
$$
\int_\O\|f(\cdot,\o)\|_{L^\be(E)}\|g(\cdot,\o)\|_{L^\be(F)}\|h(\cdot,\o)\|_{L^\be(G)}\,d\s(\o)<\be.
$$

Suppose now that $T_1$ is a bounded linear operator from $\h_2$ to $\h_1$ and $T_2$ is a bounded linear operator from $\h_3$ to $\h_2$. For a function $\Phi$ in
$L^\be(E_1)\hat\otimes_{\rm i}L^\be(E_2)\hat\otimes_{\rm i}L^\be(E_3)$ of the form \rf{ttp}, we put
\begin{align}
\label{opr}
&\int\limits_\X\int\limits_\Y\int\limits_\cZ\Phi(x,y,z)
\,d E_1(x)T_1\,dE_2(y)T_2\,dE_3(z)\\[.2cm]
\df&\int\limits_\O\left(\,\int\limits_\X f(x,\o)\,dE_1(x)\right)T_1
\left(\,\int\limits_\Y g(y,\o)\,dE_2(y)\right)T_2
\left(\,\int\limits_\cZ h(z,\o)\,dE_3(z)\right)\,d\s(\o).\nonumber
\end{align}

The following lemma from \cite{Pe8} (see also \cite{ACDS} for a different proof) shows that the definition does not depend on the choice of a representation \rf{ttp}.

\begin{lem}
\label{kor}
Suppose that $\Phi\in L^\be(E_1)\hat\otimes_{\rm i}L^\be(E_2)\hat\otimes_{\rm i}L^\be(E_3)$. Then the 
right-hand side of {\em\rf{opr}} does not depend on the choice of a representation {\em\rf{ttp}}.
\end{lem}

It is easy to see that the following inequality holds
$$
\left\|\int\limits_\X\int\limits_\Y\int\limits_\cZ\Phi(x,y,z)
\,dE_1(x)T_1\,dE_2(y)T_2\,dE_3(z)\right\|
\le\|\Phi\|_{L^\be\hat\otimes_{\rm i}L^\be\hat\otimes_{\rm i}L^\be}\cdot\|T_1\|\cdot\|T_2\|.
$$

In particular, the triple operator integral on the left-hand side of \rf{opr} can be defined if $\Phi$ belongs to the projective
tensor product $L^\be(E_1)\hat\otimes L^\be(E_2)\hat\otimes L^\be(E_3)$, i.e., $\Phi$ admits a representation
$$
\Phi(x,y,z)=\sum_{n\ge1}f_n(x)g_n(y)h_n(z),
$$
where $f_n\in L^\be(E_1)$, $g_n\in L^\be(E_2)$, $h_n\in L^\be(E_3)$ and
$$
\sum_{n\ge1}\|f_n\|_{L^\be(E_1)}\|g_n\|_{L^\be(E_2)}\|h_n\|_{L^\be(E_3)}<\be.
$$

In a similar way one can define multiple operator integrals, see \cite{Pe8}.

Recall that earlier multiple operator integrals were considered in \cite{Pa} and \cite{S}. However, in those papers the class of functions 
$\Phi$ for which the left-hand side of \rf{opr} was defined is much narrower than in the definition given above.

Multiple operator integrals arise in connection with the problem of evaluating higher order operator derivatives. It turns out that 
if $A$ is a self-adjoint operator on Hilbert space and $K$ is a bounded self-adjoint operator, then for sufficiently nice functions $\f$ on 
$\R$, the function
\bay
\label{npr}
t\mapsto \f(A+tK)
\ey
has $n$ derivatives in the norm and the $n$th derivative can be expressed in terms of multiple operator integrals.
We are going to consider this problem in \S\,8. 

Note that recently in \cite{JTT} Haagerup tensor products were used to define multiple operator integrals. However, it is
not clear whether this can lead to a broader class of functions $\f$, for which the function $f$ has an $n$th derivative and the $n$th
derivative can be expressed in terms of a multiple operator integral.

\

\section{\bf Besov spaces}
\setcounter{equation}{0}
\label{s4}

\

The purpose of this section is to give a brief introduction to the Besov spaces that play an important role in problems of perturbation theory.
We start with Besov spaces on the unit circle.

Let $1\le p,\,q\le\be$ and $s\in\R$. The Besov class $B^s_{pq}$ of functions (or
distributions) on $\T$ can be defined in the following way. Let $w$ be an infinitely differentiable function on $\R$ such
that
\bay
\label{w}
w\ge0,\quad\supp w\subset\left[\frac12,2\right],\quad\mbox{and} \quad w(x)=1-w\left(\frac x2\right)\quad\mbox{for}\quad x\in[1,2].
\ey
and $w$ is a linear function on the intervals $[1/2,1]$ and $[1,2]$.

Consider the trigonometric polynomials $W_n$, and $W_n^\#$ defined by
$$
W_n(z)=\sum_{k\in\Z}w\left(\frac{k}{2^n}\right)z^k,\quad n\ge1,\quad W_0(z)=\bar z+1+z,\quad
\mbox{and}\quad W_n^\#(z)=\ov{W_n(z)},\quad n\ge0.
$$
Then for each distribution $\f$ on $\T$,
$$
\f=\sum_{n\ge0}\f*W_n+\sum_{n\ge1}\f*W^\#_n.
$$
The Besov class $B^s_{pq}$ consists of functions (in the case $s>0$) or distributions $\f$ on $\T$
such that
\bay
\label{bes}
\big\{\|2^{ns}\f*W_n\|_{L^p}\big\}_{n\ge0}\in\ell^q\quad\mbox{and}
\quad\big\{\|2^{ns}\f*W^\#_n\|_{L^p}\big\}_{n\ge1}\in\ell^q
\ey
Besov classes admit many other descriptions. In particular, for $s>0$, the space $B^s_{pq}$ admits the
following characterization. A function $\f$ belongs to $B^s_{pq}$, $s>0$, if and only if
$$
\int_\T\frac{\|\D^n_\t \f\|_{L^p}^q}{|1-\t|^{1+sq}}d\m(\t)<\be\quad\mbox{for}\quad q<\be
$$
and
\bay
\label{pbe}
\sup_{\t\ne1}\frac{\|\D^n_\t \f\|_{L^p}}{|1-\t|^s}<\be\quad\mbox{for}\quad q=\be,
\ey
where $\m$ is normalized Lebesgue measure on $\T$, $n$ is an integer greater than $s$, and $\D_\t$ is
the difference operator: 
$$
(\D_\t \f)(\z)=\f(\t\z)-\f(\z), \quad\z\in\T.
$$

We use the notation $B_p^s$ for $B_{pp}^s$.

The spaces $\L_\a\df B_\be^\a$ form the {\it H\"older--Zygmund scale}. If $0<\a<1$, then $\f\in\L_\a$ if and only if
$$
|\f(\z)-\f(\t)|\le\const|\z-\t|^\a,\quad\z,\,\t\in\T,
$$
while $f\in\L_1$ if and only if
$$
|\f(\z\t)-2\f(\z)+\f(\z\bar\t)|\le\const|1-\t|,\quad\z,\,\t\in\T.
$$
It follows from \rf{pbe} that for $\a>0$, $\f\in\L_\a$ if and only if
$$
|(\D^n_\t\f)(\z)|\le\const|1-\t|^\a,
$$
where $n$ is a positive integer such that $n>\a$.

It is easy to see from the definition of Besov classes that the Riesz projection $\pp_+$,
$$
\pp_+\f=\sum_{n\ge0}\hat \f(n)z^n,
$$
is bounded on $B^s_{pq}$ and functions in $\big(B^s_{pq}\big)_+\df\pp_+B^s_{pq}$ admit a natural extension to analytic functions
in the unit disk $\dd$. It is well known that the functions in $\big(B^s_{pq}\big)_+$ admit the following description:
$$
\f\in \big(B^s_{pq}\big)_+\Leftrightarrow
\int_0^1(1-r)^{q(n-s)-1}\|\f^{(n)}_r\|^q_p\,dr<\be,\quad q<\be,
$$
and
$$
\f\in \big(B^s_{p\be}\big)_+\Leftrightarrow
\sup_{0<r<1}(1-r)^{n-s}\|\f^{(n)}_r\|_p<\be,
$$
where $\f_r(\z)\df \f(r\z)$ and $n$ is a nonnegative integer greater than $s$.

Besov spaces play  a significant role in many problems of operator theory and it is also important to consider Besov spaces
$B_{pq}^s$ when $p$ or $q$ can be less than $1$. Everything mentioned above also holds for arbitrary positive $p$ and $q$ provided $s>1/p-1$.

Let us proceed now to Besov spaces on the real line. We consider homogeneous Besov spaces 
$B_{pq}^s(\R)$ of functions (distributions) on $\R$.
We use the same function $w$
as in \rf{w} and define the functions $W_n$ and $W^\#_n$ on $\R$ by
$$
\F W_n(x)=w\left(\frac{x}{2^n}\right),\quad\F W^\#_n(x)=\F W_n(-x),\quad n\in\Z,
$$
where $\F$ is the {\it Fourier transform}. The Besov class $B^s_{pq}(\R)$ consists of
distributions $\f$ on $\R$ such that
$$
\{\|2^{ns}\f*W_n\|_{L^p}\}_{n\in\Z}\in\ell^q(\Z)\quad\mbox{and}
\quad\{\|2^{ns}\f*W^\#_n\|_{L^p}\}_{n\in\Z}\in\ell^q(\Z).
$$
According to this definition, the space $B^s_{pq}(\R)$ contains all polynomials. However, it is not
necessary to include all polynomials. It is natural to assume that the space $B^s_{pq}(\R)$ 
contains no polynomial of degree greater than $s-1/p$.

Besov spaces $B^s_{pq}(\R)$ admit equivalent definitions that are similar to those discussed above in the case of Besov
spaces of functions on $\T$.

We refer the reader to \cite{Pee} and \cite{Pe6} for more detailed information on Besov spaces.

\

\section{\bf Nuclearity of Hankel operators}
\setcounter{equation}{0}
\label{s5}

\

It turns out (see  \cite{Pe1}) that Hankel operators play an important role in our problems of perturbation theory.
For a function $\f$ on the unit circle $\T$, the Hankel operator $H_\f$ on the Hardy class $H^2\subset L^2$ is defined
by
$$
H_\f:H^2\to H^2_-\df L^2\ominus H^2,\quad H_\f f\df\pp_-\f f,
$$
where $\pp_-$ is the orthogonal projection onto $H^2_-$. By Nehari's theorem,
$$
\|H_\f\|=\dist_{L^\be}(\f,H^\be)
$$
(see \cite{Pe6}, Ch. 1, \S\,1).

In this paper we need the following result that describes the Hankel operator of trace class $\bS_1$.

\begin{thm}
\label{ncr}
$H_\f\in\bS_1$ if and only if $\pp_-\f\in B_1^1$.
\end{thm}

Theorem \ref{ncr} was obtained in \cite{Pe0}, see also \cite{Pe6}, Ch. 6, \S\,1.

Consider now the following class of integral operators. Let $\f\in L^\be$. The operator $\cc_\f$ on $L^2$ is defined by
$$
(\cc_\f f)(\z)=\int_\T\frac{\f(\z)-\f(\t)}{1-\bar\t\z}f(\t)\,d\m(\t).
$$

The following result can be deduced from Theorem \ref{ncr} (see \cite{Pe6}, Ch. 6, \S\.7).

\begin{thm}
\label{com}
Let $\f\in L^\be$. Then $\cc_\f\in\bS_1$ if and only if $\f\in B^1_1$.
\end{thm}

\Pf Indeed, it is easy to show that 
$$
\cc_\f f=H_\f f_+ - H^*_{\ov\f}f_-,
$$
where $f_-=\pp_-f$ and $f_+=f-f_-$. Theorem \ref{com} follows now immediately from Theorem \ref{ncr}. $\bl$

\

\section{\bf Operator Lipschitz and operator differentiable functions. \\ Sufficient conditions}
\setcounter{equation}{0}
\label{s6}

\

In this section we discuss  sufficient conditions for a function on the unit circle or on the real line to be 
operator Lipschitz or operator differentiable. 
We begin with unitary operators. 

The following lemma gives us an estimate for the norm of $\dg\f$ in the projective tensor product
$L^\be\hat\otimes L^\be$ in the case of trigonometric polynomials  $\f$. It was obtained in \cite{Pe1}.
We give here a slightly modified proof given in \cite{Pe8}.

\begin{lem}
\label{cl}
Let $\f$ be a trigonometric polynomial of degree $m$. Then
\bay
\label{cru}
\|\dg\f\|_{L^\be\hat\otimes L^\be}\le\const m\|\f\|_{L^\be}.
\ey
\end{lem}

\Pf First of all, it is evident that it suffices to consider the case when $m=2^l$. 
Next, it suffices to prove the result for analytic polynomials $\f$ (i.e., linear combinations of $z^j$ with $j\ge0$). 
Indeed, if \rf{cru} holds for analytic polynomials, then it obviously also holds for conjugate trigonometric polynomials. Let now $\f$
be an arbitrary trigonometric polynomial of degree $2^l$. We have
$$
\f=\sum_{j=1}^l\f*W^\#_j+\sum_{j=0}^l\f*W_j
$$
(see \S\,4). 
Applying \rf{cru} to each term of this expansion, we obtain
\begin{align*}
\|\dg\f\|_{L^\be\hat\otimes L^\be}&\le\sum_{j=1}^l\big\|\dg\big( \f*W_j^\#\big)\big\|_{L^\be\hat\otimes L^\be}
+\sum_{j=0}^l\|\dg( \f*W_j)\|_{L^\be\hat\otimes L^\be}\\[.2cm]
&\le\const\left(\sum_{j=1}^l2^j\big\|\f*W_j^\#\big\|_{L^\be}+\sum_{j=0}^l2^j\|\f*W_j\|_{L^\be}\right)\\[.2cm]
&\le\const\sum_{j=0}^l2^j\|\f\|_{L^\be}\le\const2^l\|\f\|_{L^\be}.
\end{align*}

Assume now that $\f$ is an analytic polynomial of degree $m$. It is easy to see that
$$
(\dg\f)(z_1,z_2)=\sum_{j,k\ge0}\hat\f(j+k+1)z_1^jz_2^k.
$$
We have
$$
\sum_{j,k\ge0}\hat\f(j+k+1)z_1^jz_2^k=\sum_{j,k\ge0}\a_{jk}\hat\f(j+k+1)z_1^jz_2^k+\sum_{j,k\ge0}\b_{jk}\hat\f(j+k+1)z_1^jz_2^k,
$$
where
$$
\a_{jk}=\left\{\begin{array}{ll}\frac12,&j=k=0,\\[.2cm]
\frac{j}{j+k},& j+k\ne0
\end{array}\right.
$$
and
$$
\b_{jk}=\left\{\begin{array}{ll}\frac12,&j=k=0,\\[.2cm]
\frac{k}{j+k},& j+k\ne0.
\end{array}\right.
$$
Clearly, it is sufficient to estimate
$$
\left\|\sum_{j,k\ge0}\a_{jk}\hat\f(j+k+1)z_1^jz_2^k\right\|_{L^\be\hat\otimes L^\be}.
$$
It is easy to see that
$$
\sum_{j,k\ge0}\a_{jk}\hat\f(j+k+1)z_1^jz_2^k=
\sum_{k\ge0}\left(\Big(\big((S^*)^{k+1}\f\big)*\sum_{j\ge0}\a_{jk}z^j\Big)(z_1)\right)z_2^k,
$$
where $S^*$ is backward shift, i.e., $(S^*)^k\f=\pp_+\bar z^k\f$. 

Thus
$$
\left\|\sum_{j,k\ge0}\a_{jk}\hat\f(j+k+1)z_1^jz_2^k\right\|_{L^\be\hat\otimes L^\be}\le
\sum_{k\ge0}\left\|\big((S^*)^{k+1}\f\big)*\sum_{j\ge0}\a_{jk}z^j\right\|_{L^\be}.
$$
Put
$$
Q_k(z)=\sum_{i\ge k}\frac{i-k}{i}z^i,\quad k>0,\quad \mbox{and}\quad
Q_0(z)=\frac12+\sum_{i\ge1}z^i.
$$
Then it is easy to see that
$$
\left\|\big((S^*)^{k+1}\f\big)*\sum_{j\ge0}\a_{jk}z^j\right\|_{L^\be}=\|\psi*Q_{k}\|_{L^\be},
$$
where $\psi=S^*\f$,
and so 
$$
\left\|\sum_{j,k\ge0}\a_{jk}\hat\f(j+k+1)z_1^jz_2^k\right\|_{L^\be\hat\otimes L^\be}
\le\sum_{k\ge0}\|\psi*Q_{k}\|_{L^\be}.
$$
Consider the function $r$ on $\R$ defined by
$$
r(x)=\left\{\begin{array}{ll}1,&|x|\le1,\\[.2cm]\frac1{|x|},&|x|\ge1.
\end{array}\right.
$$
It is easy to see that the Fourier transform $\F r$ of $r$ belongs to $L^1(\R)$.
Define the functions $R_n$, $n\ge1$, on $\T$ by
$$
R_k(\z)=\sum_{j\in\Z}r\left(\frac jk\right)\z^j.
$$
An elementary estimate obtained in Lemma 4.3 of  \cite{Pe8} shows that
$$
\|R_k\|_{L^1}\le\const.
$$
It is easy to see that
for $f\in H^\be$, we have
$$
\|f*Q_k\|_{L^\be}=\|f-f*R_k\|_{L^\be}\le\|f\|_{L^\be}+\|f*R_k\|_{L^\be}\le\const\|f\|_{L^\be}.
$$
Thus
$$
\sum_{k\ge0}\|\psi*Q_k\|_{L^\be}=\sum_{k=0}^m\|\psi*Q_k\|_{L^\be}\le\const m\|\psi\|_{L^\be}\le\const m\|\f\|_{L^\be}.\quad\bl
$$

The following result was obtained in \cite{Pe1}.

\begin{thm}
\label{du}
Let $\f\in B^1_{\be1}$. Then $\dg\f\in C(\T)\hat\otimes C(\T)$ and
$$
\|\dg\f\|_{L^\be\hat\otimes L^\be}\le\const\|\f\|_{B^1_{\be1}}.
$$
\end{thm}

\Pf We have 
$$
\f=\sum_{j>0}\f*W_j^\#+\sum_{j\ge0}\f*W_j.
$$
By Lemma \ref{cru}, each of the functions $\dg\big(\f*W_j^\#\big)$  and $\dg(\f*W_j)$ belongs to $C(\T)\hat\otimes C(\T)$ and
\begin{align*}
\sum_{j>0}\big\|\dg\big(\f*W_j^\#\big)\big\|_{L^\be\hat\otimes L^\be}&+\sum_{j\ge0}\big\|\dg(\f*W_j)\big\|_{L^\be\hat\otimes L^\be}\\[.2cm]
&\le\const\left(\sum_{j>0}2^j\big\|\f*W_j^\#\big\|_{L^\be}+\sum_{j\ge0}2^j\|\f*W_j\|_{L^\be}\right)\\[,2cm]
&\le\const\|\f\|_{B^1_{\be1}}.\qquad\bl
\end{align*}

It follows from Theorem  \ref{du} that for $\f\in B^1_{\be1}$, the divided difference 
$\dg\f$ belongs to the space  $\fM(E,F)$ of Schur multipliers with respect to arbitrary Borel spectral measures $E$ and $F$ on $\T$
(see \S\,2). By the Birman--Solomyak formula for unitary operators, we have
\bay
\label{BSu}
\f(U)-\f(V)=\iint\limits_{\T\times\T}\frac{\f(\z)-\f(\t)}{\z-\t}\,dE_U(\z)\,(U-V)\,dE_V(\t),
\ey
which implies the following result:

\begin{thm}
\label{uoL}
Let $\f\in B^1_{\be1}$. Then $\f$ is operator Lipschitz, i.e.,
$$
\|\f(U)-\f(V)\|\le\const\|U-V\|,
$$
for unitary operators $U$ and $V$ on Hilbert space.
\end{thm}

\Pf It follows from \rf{BSu} that
\begin{align*}
\|\f(U)-\f(V)\|&\le\|\dg\f\|_{\fM(E_U,E_V)}\|U-V\|\\[.2cm]
&\le\|\dg\f\|_{L^\be\hat\otimes L^\be}\|U-V\|
\le\const\|\f\|_{B^1_{\be1}}\|U-V\|.\quad\bl
\end{align*}

Let us now show that the condition $\f\in B^1_{\be1}$ also implies that $\f$ is operator differentiable. 

\begin{thm}
\label{odu}
Let $\f$ be a function on $\T$ of class $B^1_{\be1}$. If $A$ is a bounded self-adjoint operator and $U$ is a unitary operator,
and $U_s\df e^{{\rm i}sA}U$, then
the function 
\bay
\label{sdi}
s\mapsto\f(U_s)
\ey
 is differentiable in the norm 
and
\bay
\label{DKu}
\frac{d}{ds}\big(\f(U_s)\big)\Big|_{s=o}
={\rm i}\left(\iint\frac{\f(\z)-\f(\t)}{\z-\t}\,dE_U(\z)A\,dE_U(\t)\right)U.
\ey
Moreover, the map
\bay
\label{Fdm}
A\mapsto \f\big(e^{{\rm i}A}U\big)
\ey
defined on the space of bounded self-adjoint operators is differentiable in the sense of Fr\'echet.
\end{thm}

\Pf Let us prove that the function \rf{sdi} is norm differentiable and that formula
\rf{DKu} holds.
By Theorem \ref{du}, there exist continuous functions $f_n$ and $g_n$ on $\T$ such that
$$
(\dg\f)(\z,\t)=\sum_{n\ge1}f_n(\z)g_n(\t),\quad\z,\,\t\in\T,
$$
and
\bay
\label{tni}
\sum_{n\ge1}\|f_n\|_{L^\be}\|g_n\|_{L^\be}<\be.
\ey
By the Birman--Solomyak formula \rf{BSu},
\begin{align*}
\f(U_s)-\f(U)&=\iint\limits_{\T\times\T}(\dg\f)(\z,\t)\,dE_{U_s}(\z)(U_s-U)\,dE_U(\t)\\[.2cm]
&=\sum_{n\ge1}f_n(U_s)(U_s-U)g_n(U).
\end{align*}
On the other hand,
$$
{\rm i}\left(\iint\frac{\f(\z)-\f(\t)}{\z-\t}\,dE_U(\z)A\,dE_U(\t)\right)U=
{\rm i}\sum_{n\ge1}f_n(U)Ag_n(U)U.
$$
We have 
\begin{align*}
\frac1s\big(\f(U_s)-\f(U)\big)&-{\rm i}\sum_{n\ge1}f_n(U)Ag_n(U)U\\[.2cm]
&=
\sum_{n\ge1}\left(\frac1s f_n(U_s)(U_s-U)g_n(U)-{\rm i}f_n(U)Ag_n(U)U\right)\\[.2cm]
&=\sum_{n\ge1}\left(\frac1s f_n(U_s)(U_s-U)g_n(U)-\frac1sf_n(U)(U_s-U)g_n(U)\right)\\[.2cm]
&+\sum_{n\ge1}\left(\frac1s f_n(U)(U_s-U)g_n(U)-{\rm i}f_n(U)Ag_n(U)U\right).
\end{align*}
Clearly,
$$
\left\|\frac1s(U_s-U)\right\|\le\const.
$$
Since, $f_n\in C(\T)$, it is easy to see that 
$$
\lim_{s\to0}\|f_n(U_s)-f_n(U)\|=0.
$$
It follows now easily from \rf{tni} that
$$
\lim_{s\to0}\left\|
\sum_{n\ge1}\left(\frac1s f_n(U_s)(U_s-U)g_n(U)-\frac1sf_n(U)(U_s-U)g_n(U)\right)
\right\|=0.
$$
On the other hand, it is easy to see that 
$$
\lim_{s\to\0}\left\|
\frac1s(U_s-U)-{\rm i}AU
\right\|=0
$$
and again, it follows from \rf{tni} that
$$
\lim_{s\to0}\left\|
\sum_{n\ge1}\left(\frac1s f_n(U)(U_s-U)g_n(U)-{\rm i}f_n(U)Ag_n(U)U\right)
\right\|=0
$$
which proves that the function \rf{sdi} is norm differentiable and \rf{DKu} holds. 

One can easily see that the same reasoning also shows that the map \rf{Fdm} is differentiable in the sense of Fr\'echet. $\bl$

The above results of this section were obtained in \cite{Pe1}. Recall that earlier Birman and Solomyak proved in \cite{BS3} the same results for functions $\f$ whose derivatives belong to the H\"older class $\L_\a$ with some $\a>0$.

In the case of differentiability in the Hilbert--Schmidt norm, the following result was proved by Birman and Solomyak in \cite{BS3}.

\begin{thm}
\label{BSHS}
Let $\f\in C^1(\T)$. If under the hypotheses of Theorem {\em\ref{odu}} the self-adjoint operator $A$ belongs to the Hilbert--Schmidt class
$\bS_2$, then formula {\em\rf{DKu}} holds in the Hilbert--Schmidt norm.
\end{thm}

Let us state now similar results for (not necessarily bounded) self-adjoint operators. The following result shows that functions 
in $B_{\be1}^1(\R)$ are operator Lipschitz.

\begin{thm}
\label{oLsa}
Let $\f$ be a function on $\R$ of class $B_{\be1}^1(\R)$ and let $A$ and $B$ be self-adjoint operators such that
$A-B$ is bounded. Then the operator $\f(A)-\f(B)$ is bounded and
$$
\|\f(A)-\f(B)\|\le\const\|\f\|_{B^1_{\be1}(\R)}\|A-B\|.
$$
\end{thm}

\begin{thm}
\label{odsa}
Let $\f\in B_{\be1}^1(\R)$. Suppose that $A$ is a self-adjoint operator (not necessarily bounded) and $K$ is a bounded self-adjoint
operator. Then the function
$$
t\mapsto f(A+tK)-f(A)
$$
is norm differentiable and
$$
\frac{d}{dt}f(A+tK)\Big|_{t=0}=\iint\limits_{\R\times\R}\frac{\f(x)-\f(y)}{x-y}\,dE_A(x)K\,dE_A(y).
$$
Moreover, the map
$$
K\mapsto f(A+K)-f(A)
$$
defined on the space of bounded self-adjoint operators is differentiable in the sense of Fr\'echet.
\end{thm}

We refer the reader to \cite{Pe4} and \cite{Pe8} for the proofs of Theorems \ref{oLsa} and \ref{odsa}.

\

\section{\bf Operator Lipschitz and operator differentiable functions. \\ Necessary conditions}
\setcounter{equation}{0}
\label{s7}

\

\begin{thm}
\label{ncL}
Let $\f$ be a continuously differentiable function on $\T$. If $\f$ is operator Lipschitz, then $\f\in B_1^1$.
\end{thm}

Note that the condition $\f\in B_1^1$ implies that
$$
\sum_{n\ge0}|\widehat{\f'}(2^n)|<\be.
$$
This follows easily from \rf{bes}. On the other hand, it is well known that for an arbitrary sequence $\{c_n\}_{n\ge0}$ in $\ell^2$,
there exists $\f\in C^1(\T)$ such that 
$$
\widehat{\f'}(2^n)=c_n,\quad n\ge0.
$$
Thus the condition $\f\in C^1(\T)$ is not sufficient for $\f$ to be operator Lipschitz.

\medskip

{\bf Proof.} Let $U$ be multiplication by $z$ on $L^2$ (with respect to Lebesgue measure) and let $A$ be a self-adjoint operator on 
$L^2$ of class $\bS_2$. Put $V_t=e^{{\rm i}tA}U$, $t\in\R$. It is easy to see that
$$
\frac1t\|V_t-U\|\le\const\|A\|,
$$
and since $\f$ is operator Lipschitz, we have
$$
\left\|\frac1t\big(\f(V_t)-\f(U)\big)\right\|\le\const\|A\|.
$$
By Theorem \ref{BSHS},
$$
\lim_{t\to0}\frac1t\big(\f(V_t)-\f(U)\big)={\rm i}\left(\iint\frac{\f(\z)-\f(\t)}{\z-\t}\,dE_U(\z)A\,dE_U(\t)\right)U
$$
in the Hilbert--Schmidt norm. It follows that
$$
\left\|\iint\frac{\f(\z)-\f(\t)}{\z-\t}\,dE_U(\z)A\,dE_U(\t)\right\|\le\const\|A\|.
$$
This means that the divided difference $\dg\f$ is a Schur multiplier in $\fM(E_U,E_U)$.

Consider now the class $\frak N$ of kernel functions of trace class integral operators on $L^2$ with respect to Lebesgue measure.
By Theorem \ref{zam},
$$
k\in{\frak N}\imp(\dg\f) k\in{\frak N}.
$$
Put now
$$
k(\z,\t)\df\t.
$$
Clearly, the integral operator with kernel function $k$ is a rank one operator.
We have
$$
\big((\dg\f) k\big)(\z,\t)=\frac{\f(\z)-\f(\t)}{1-\bar\t\z},\quad\z,\,\t\in\T.
$$
Thus $\cc_\f\in\bS_1$ and it follows now from Theorem \ref{com} that $\f\in B_1^1$. $\bl$

\medskip

{\bf Remark.} It is easy to see that the reasoning given in the proof of Theorem \ref{ncL} also gives the following result:

{\it Suppose that $\f\in C^1(\T)$ and the divided difference $\dg\f$ is not a Schur multiplier in $\fM(E_U,E_U)$
(or, in other words, $\dg\f$ is not a multiplier of the class $\frak N$ of kernel functions of trace class integral operators on 
$L^2(\m)$). Then $\f$ is not operator Lipschitz.}

\medskip

Theorem \ref{ncL} was proved in \cite{Pe1}. A more elaborate application of the nuclearity criterion for Hankel operators allowed the author to obtain in \cite{Pe1} a stronger necessary condition. To state it, we introduce the class $\cL$.

\medskip

{\bf Definition.} A bounded function $\f$ on $\T$ is said to belong to $\cL$ if the Hankel operators  $H_\f$ and $H_{\bar\f}$ map
the Hardy class $H^1$ into the Besov space $B_1^1$, i.e.,
$$
\pp_-\f g\in B_1^1\quad\mbox{and}\quad \pp_-\bar\f f\in B_1^1,
$$
whenever $f\in H^1$.

\medskip

It is easy to see that $\cL\subset B^1_1$. 

The following result was obtained in \cite{Pe1}.

\begin{thm}
\label{cL}
Let $\f$ be an operator Lipschitz function of class $C^1(\T)$. Then $\f\in\cL$.
\end{thm}

The proof of Theorem \ref{cL} is given in \cite{Pe1}. It is based on the nuclearity criterion for Hankel operators, see Theorem \ref{ncr}.
Actually, it is shown in \cite{Pe1} that if $\f\in C^1(\T)\setminus\cL$, then $\dg\f$ is not a multiplier of the class $\frak N$.

S. Semmes observed (see the proof in \cite{Pe5}) that $\f\in\cL$ if and only if the measure
$$
\|{\rm Hess}\,\f\|\,dx\,dy
$$
is a Carleson measure in the unit disk $\dd$, where ${\rm Hess}\,\f$ is the Hessian matrix of the harmonic
extension of $\f$ to the unit disk.

M. Frazier observed that actually $\cL$ is the Triebel--Lizorkin space $F^1_{\be1}$. Note that the definition of the
Triebel--Lizorkin spaces $\dot F^s_{pq}$ on $\R^n$ for $p=\be$ and  $q>1$ can be found in \cite{T}, \S\,5.1.
A definition for all $q>0$, which is equivalent to Triebel's
definition when $q>1$, was given by Frazier and Jawerth in \cite{FrJ}. 
Their approach did not use harmonic extensions, but a
straightforward exercise in comparing kernels shows that Frazier and
Jawerth's definition of $\dot F^1_{\infty1}$ is equivalent to the
definition requiring $\|{\rm Hess\,\f}\| dx dy$ to be a Carleson measure on the
upper half-space.  Our space $\cL$ is the analogue for the unit disc.

The condition $\f\in\cL$ (and a fortiori the condition $\f\in B_1^1$) is also a necessary condition for the function $\f$ to
be operator differentiable.

Similar results also hold in the case of functions of self-adjoint operators:

\begin{thm}
\label{oLs}
Let $\f$ be a continuously differentiable function on $\R$. If $\f$ is locally operator Lipschitz, then
$\f$ belongs to $B_1^1(\R)$ locally.
\end{thm}

Note that the latter property means that the restriction of $\f$ to any finite interval coincides with the restriction to this interval of a function of class $B_1^1(\R)$.

\begin{thm}
\label{uoLs}
Let $\f$ be a continuously differentiable function on $\R$. If $\f$ is operator Lipschitz, then
$\f$ belongs to the class $\cL(\R)$.
\end{thm}

Note that the class $\cL(\R)$ can be defined by analogy with the class $\cL$ of functions on $\T$. The same description in terms of Carleson measures also holds. Theorem \ref{oLs} can be improved: if $\f$ is locally operator Lipschitz, then $\f$ must belong to $\cL(\R)$ locally.

Theorems \ref{oLs} and \ref{uoLs} were proved in \cite{Pe1} and \cite{Pe4}. 

To conclude this section, we mention that the same necessary conditions also hold for operator differentiability.

\

\section{\bf Higher order operator derivatives}
\setcounter{equation}{0}
\label{s8}

\

For a function $\f$ on the circle the  {\it divided differences $\dg^k\f$ of order $k$} are defined inductively as follows:
$$
\dg^0\f\df\f;
$$
if $k\ge1$, then in the case when $\l_1,\l_2,\cdots,\l_{k+1}$ are distinct points in $\T$,
$$
(\dg^{k}\f)(\l_1,\cdots,\l_{k+1})\df
\frac{(\dg^{k-1}\f)(\l_1,\cdots,\l_{k-1},\l_k)-
(\dg^{k-1}\f)(\l_1,\cdots,\l_{k-1},\l_{k+1})}{\l_{k}-\l_{k+1}}
$$
(the definition does not depend on the order of the variables). Clearly,
$$
\dg\f=\dg^1\f.
$$
If $\f\in C^k(\T)$, then $\dg^{k}\f$ extends by continuity to a function defined for all points $\l_1,\l_2,\cdots,\l_{k+1}$.

The following result was established in  \cite{Pe8}.

\begin{thm}
\label{ndd}
Let $n\ge1$ and let $\f$ be a function on $\T$ of class $B_{\be1}^n$. Then 
$\dg^n\f$ belongs to the projective tensor product $\underbrace{C(\T)\hat\otimes \cdots\hat\otimes C(\T)}_{n+1}$
and 
\bay
\label{nrr}
\|\dg^n\f\|_{L^\be\hat\otimes \cdots\hat\otimes L^\be}\le\const\|\f\|_{B_{\be1}^n}.
\ey
\end{thm}

The constant on the right-hand side of \rf{nrr} can depend on $n$.

As in the case of double operator integrals, the following lemma gives us a crucial estimate.

\begin{lem}
\label{rrp}
Let $n$ and $m$ be a positive integers and let $\f$ be a trigonometric polynomial of degree $m$.
Then
\bay
\label{Bin}
\|\dg^n\f\|_{L^\be\hat\otimes \cdots\hat\otimes L^\be}\le\const m^n\|\f\|_{L^\be}.
\ey
\end{lem}

Note that the constant on the right-hand side of \rf{Bin} can depend on $n$, but does not depend on $m$. 

The proof of Lemma \ref{rrp} is based on the same ideas as the proof of Lemma \ref{cl}. We refer the reader to \cite{Pe8} for the proof of Lemma \ref{rrp}.

We deduce now Theorem \ref{ndd} from Lemma \ref{rrp}.

\medskip

{\bf Proof of Theorem \ref{ndd}.} We have 
$$
\f=\sum_{j>0}\f*W_j^\#+\sum_{j\ge0}\f*W_j.
$$
By Lemma \ref{rrp}, 
\begin{align*}
\|\dg^n\f\|_{L^\be\hat\otimes \cdots\hat\otimes L^\be}&\le
\sum_{j>0}\big\|\dg^n(\f*W_j^\#)\big\|_{L^\be\hat\otimes\cdots\hat\otimes L^\be}+
\sum_{j\ge0}\big\|\dg^n(\f*W_j)\big\|_{L^\be\hat\otimes\cdots\hat\otimes L^\be}\\[.2cm]
&\le\const\left(\sum_{j>0}2^{jn}\|\f*W_j^\#\|_{L^\be}+\sum_{j\ge0}2^{jn}\|\f*W_j\|_{L^\be}\right)\\[,2cm]
&\le\const\|\f\|_{B^n_{\be1}}.\qquad\bl
\end{align*}

Suppose now that $U$ is a unitary operator and $A$ is a bounded self-adjoint operator on Hilbert space. Consider the family of operators
$$
U_t=e^{{\rm i}tA}U.\quad t\in\R.
$$
The following result was proved in \cite{Pe8}.

\begin{thm}
\label{npu}
Let $\f\in B^n_{\be1}$. Then the function
$$
t\mapsto\f(U_t)
$$
has $n$ derivatives in the norm and
\begin{align*}
&\frac{d^n}{dt^n}\big(\f(U_t)\big)\Big|_{s=0}\\[.2cm]
=&{\rm i}^nn!\left(\underbrace{\int\cdots\int}_{n+1}(\dg^n\f)(\z_1,\cdots,\z_{n+1})
\,dE_U(\z_1)A\cdots A\,dE_U(\z_{n+1})\right)U^n.
\end{align*}
\end{thm}

Similar results hold for self-adjoint operators. The following results can be found in \cite{Pe8}:

\begin{thm}
\label{rrn}
Let $\f\in B^n_{\be1}(\R)$. Then $\dg^n\f$ belongs to the integral projective tensor product 
$\underbrace{L^\be\hat\otimes \cdots\hat\otimes L^\be}_{n+1}$.
\end{thm}

\begin{thm}
\label{nprs}
Suppose that $\f\in B^n_{\be1}(\R)\cap B^1_{\be1}(\R)$. Let $A$ be a self-adjoint operator and let $K$ be a bounded self-adjoint
operator on Hilbert space. Then the function
\bay
\label{fAK}
t\mapsto\f(A+tK)
\ey
has an $n$th derivative in the norm and
\bay
\label{nopr}
\frac{d^n}{dt^n}\big(\f(A+tK)\big)\Big|_{s=0}
\!=\!n!\!\underbrace{\int\cdots\int}_{n+1}(\dg^n\f)(x_1,\cdots,x_{n+1})
\,dE_U(x_1)A\cdots A\,dE_U(x_{n+1}).
\ey
\end{thm}

Note that under the hypotheses of Theorem \ref{nprs}, the function \rf{fAK} has all the derivatives in the norm up to order $n$.
However, in the case of unbounded self-adjoint operators it can happen that the $n$th derivative exists in the norm, but lower order derivatives do not exist in the norm. For example, if $\f\in B^2_{\be1}(\R)$, but $\f\not\in B^1_{\be1}(\R)$, then it can happen that
the function \rf{fAK} does not have the first derivative in the norm, but it is possible to interpret its second derivative so that the second derivative exists in the norm and can be computed by formula \rf{nopr}; see detailed comments in \cite{Pe8}.

Earlier sufficient conditions for the function \rf{fAK} to have $n$ derivatives in the norm and satisfy \rf{nopr} were found
in \cite{S}. However, the conditions found in \cite{S} were much more restrictive.

\

\section{\bf The case of contractions}
\setcounter{equation}{0}
\label{s9}

\

Let $T$ be a contraction (i.e., $\|T\|\le1$) on Hilbert space. 
Von Neumann's inequality (see \cite{SNF}) says that for an arbitrary analytic polynomial $\f$,
$$
\|\f(T)\|\le\max_{|\z|\le1}|\f(\z)|.
$$
This allows one to define the functional calculus
$$
\f\mapsto\f(T)
$$
on the disk-algebra $C_A$. 

In this case we consider the questions of the behavior of $\f(T)$ under perturbations of $T$. As in the case of
unitary operators a function $\f\in C_A$ is called {\it operator Lipschitz} if
$$
\|\f(T)-\f(R)\|\le\const\|T-R\|
$$
for arbitrary contractions $T$ and $R$. We also consider differentiability properties. 

For contractions $T$ and $R$, we consider the one-parameter family of contractions
$$
T_t=(1-t)T+tR,\quad\0\le t\le1,
$$
and we study differentiability properties of the map
\bay
\label{opf}
t\mapsto\f(T_t)
\ey
for a given function $\f$ in $C_A$.

It was observed in \cite{Pe2} that if $\f$ is an analytic function in $\big(B_{\be1}^1\big)_+$, then $\f$ is operator Lipschitz.
To prove this, double operator integrals with respect to semi-spectral measures were used.

Recently in \cite{KS} it was proved that if $\f\in C_A$, the the following are equivalent:

(i) $\|\f(U)-\f(V)\|\le\const\|U-V\|$ for arbitrary unitary operators $U$ and $V$;

(ii) $\|\f(T)-\f(R)\|\le\const\|T-R\|$ for arbitrary contractions $T$ and $R$.

In \cite{Pe9} it was shown that the same condition $\big(B_{\be1}^1\big)_+$ implies that the function \rf{opf}
is differentiable in the norm and the derivative can be expressed in terms of double operator integrals with respect to semi-spectral measures. It was also established in \cite{Pe9} that under the condition $\f\in B^n_{\be1}$, the function \rf{opf} is $n$ times differentiable in the norm and the $n$th derivative can be expressed in terms of a multiple operator integral with respect to semi-spectral measures.

\medskip

{\bf Definition.} Let $\h$ be a Hilbert space and let $(\X,\B)$ be a measurable space.
A map $\E$ from $\B$ to the algebra $B(\h)$ of all bounded operators on $\h$ is called a {\it semi-spectral measure}
if 
$$
\E(\D)\ge\0,\quad\D\in\B,
$$
$$
\E(\varnothing)=\0\quad\mbox{and}\quad\E(\X)=I,
$$
and for a sequence $\{\D_j\}_{j\ge1}$ of disjoint sets in $\B$,
$$
\E\left(\bigcup_{j=1}^\be\D_j\right)=\lim_{N\to\be}\sum_{j=1}^N\E(\D_j)\quad\mbox{in the weak operator topology}.
$$

\medskip

If $\K$ is a Hilbert space, $(\X,\B)$ is a measurable space,  $E:\B\to B(\K)$ is a spectral measure, and $\h$ is
a subspace of $\K$, then it is easy to see that the map $\E:\B\to B(\h)$ defined by
\bay
\label{dil}
\E(\D)=P_\h E(\D)\big|\h,\quad\D\in\B,
\ey
is a semi-spectral measure. Here $P_\h$ stands for the orthogonal projection onto $\h$.

Naimark proved in \cite{Na}  (see also \cite{SNF}) that all semi-spectral measures can be obtained in this way, i.e.,
a semi-spectral measure is always a {\it compression} of a spectral measure. A spectral measure $E$ satisfying \rf{dil} is called a {\it spectral dilation of the semi-spectral measure} $\E$.

A spectral dilation $E$ of a semi-spectral measure $\E$ is called {\it minimal} if 
$$
\K=\clos\spn\{E(\D)\h:~\D\in\B\}.
$$

It was shown in \cite{MM} that if $E$ is a minimal spectral dilation of a semi-spectral measure $\E$, then
$E$ and $\E$ are mutually absolutely continuous and all minimal spectral dilations of a semi-spectral measure are isomorphic in the natural sense.

If $\f$ is a bounded complex-valued measurable function on $\X$ and $\E:\B\to B(\h)$ is a semi-spectral measure, then the integral
\bay
\label{iss}
\int_\X f(x)\,d\E(x)
\ey
can be defined as
\bay
\label{voi}
\int_\X f(x)\,d\E(x)=\left.P_\h\left(\int_\X f(x)\,d E(x)\right)\right|\h,
\ey
where $E$ is a spectral dilation of $\E$. It is easy to see that the right-hand side of \rf{voi} does not depend on the choice
of a spectral dilation. The integral \rf{iss} can also be computed as the limit of sums
$$
\sum f(x_\a)\E(\D_\a),\quad x_\a\in\D_\a,
$$
over all finite measurable partitions $\{\D_\a\}_\a$ of $\X$.

If $T$ is a contraction on a Hilbert space $\h$, then by the Sz.-Nagy dilation theorem
(see \cite{SNF}),  $T$ has a unitary dilation, i.e., there exist a Hilbert space $\K$ such that
$\h\subset\K$ and a unitary operator $U$ on $\K$ such that
\bay
\label{DT}
T^n=P_\h U^n\big|\h,\quad n\ge0,
\ey
where $P_\h$ is the orthogonal projection onto $\h$. Let $E_U$ be the spectral measure of $U$.
Consider the operator set function $\E$ defined on the Borel subsets of the unit circle $\T$ by
$$
\E(\D)=P_\h E_U(\D)\big|\h,\quad\D\subset\T.
$$
Then $\E$ is a semi-spectral measure. It follows immediately from
\rf{DT} that 
\bay
\label{step}
T^n=\int_\T \z^n\,d\E(\z)=P_\h\int_\T\z^n\,dE_U(\z)\Big|\h,\quad n\ge0.
\ey
Such a semi-spectral measure $\E$ is called a {\it semi-spectral measure} of $\T$.
Note that it is not unique. To have uniqueness, we can consider a minimal unitary dilation $U$ of $T$,
which is unique up to an isomorphism (see \cite{SNF}).

It follows easily from  \rf{step} that 
$$
\f(T)=P_\h\int_\T\f(\z)\,dE_U(\z)\Big|\h
$$
for an arbitrary function $\f$ in the disk-algebra $C_A$.

In \cite{Pe9} double operator integrals and multiple operator integrals with respect to semi-spectral measures were introduced.

Suppose that $(\X_1,\B_1)$ and $(\X_2,\B_2)$ are measurable spaces, and
$\E_1:\B_1\to B(\h_1)$ and $\E_2:\B_2\to B(\h_2)$ are semi-spectral measures.
Then double operator integrals
$$
\iint\limits_{\X_1\times\X_2}\Phi(x_1,x_2)\,d\E_1(x_1)Q\,d\E_2(X_2).
$$
were defined in \cite{Pe9} in the case when $Q\in\bS_2$ and $\Phi$ is a bounded measurable function and in the case when $Q$ is a bounded linear operator and $\Phi$ belongs to the integral projective tensor product of the spaces $L^\be(\E_1)$
and $L^\be(\E_2)$.

Similarly, multiple operator integrals with respect to semi-spectral measures were defined in \cite{Pe9} for functions that belong to the integral projective tensor product of the corresponding $L^\be$ spaces.

Let us now state the results obtained in \cite{Pe9}.

For a contraction $T$ on Hilbert space, we denote by $\E_T$ a semi-spectral measure of $\T$.
Recall that if $\f'\in C_A$, the function $\dg\f$ extends to the diagonal
$$
\bs{\D}\df\big\{(\z,\z):~\z\in\T\big\}
$$
by continuity: $(\dg\f)(\z,\z)=\f'(\z)$, $\z\in\T$.

\begin{thm}
\label{BSTc}
Let $\f\in\left(B_{\be1}^1\right)_+$. Then for contractions $T$ and $R$ on Hilbert space the following formula holds:
\bay
\label{BSFC}
\f(T)-\f(R)=\iint\limits_{\T\times\T}\frac{\f(\z)-\f(\t)}{\z-\t}\,d\E_T(\z)\,(T-R)\,d\E_R(\t).
\ey
\end{thm}

\begin{thm}
\label{S2}
Let $\f$ be a function analytic in $\dd$ such that $\f'\in C_A$. 
If $T$ and $R$ are contractions such that $T-R\in\bS_2$, then formula
{\em\rf{BSFC}} holds.
\end{thm}

{\bf Remark.} Theorem \ref{S2} can be extended to the more general case when
 $\f'\in H^\be$. In this case we can define $\dg\f$ to be zero on the diagonal $\bs{\D}$.

\medskip

The following result is an immediate consequence of the above remark; it was obtained recently in \cite{KS} by a
completely different method.

\begin{cor}
\label{hsp}
Suppose that $\f$ is a function analytic in $\dd$ such that $\f'\in H^\be$. If $T$ and $R$ are contractions on
Hilbert space such that $T-R\in\bS_2$, then
$$
\f(T)-\f(R)\in\bS_2
$$
and
$$
\|\f(R)-\f(T)\|_{\bS_2}\le\|\f'\|_{H^\be}\|T-R\|_{\bS_2}.
$$
\end{cor}

We proceed now to the differentiability problem.
Let $T$ and $R$ be contractions on Hilbert space and let $\f\in C_A$. We are interested in differentiability properties of the function 
\rf{opf}.

Let $\E_t$ be a semi-spectral measure of $T_t$ on the unit circle $\T$,
i.e.,
$$
T_t^n=\int\limits_\T\z^n\,d\E_t(\z),\quad n\ge0.
$$
Put $\E\df\E_0$.

The following results were established in \cite{Pe9}.

\begin{thm}
\label{per}
Suppose that $\f\in\big(B^1_{\be1}\big)_+$. Then the function {\em\rf{opf}} is differentiable 
in the norm and
$$
\frac{d}{ds}\f(T_s)\Big|_{s=t}=\iint\limits_{\T\times\T}\frac{\f(\z)-\f(\t)}{\z-\t}\,d\E_t(\z)\,(R-T)\,d\E_t(\t).
$$
\end{thm}

Note that the same result holds in the case when $T-R\in\bS_2$ and $\f$ is an analytic function in $\dd$ such that $\f'\in C_A$. In this case the derivative exists in the Hilbert --Schmidt norm; see \cite{Pe9}.

We conclude this section with an analog of Theorem \ref{nprs} for contractions.

\begin{thm}
\label{n}
Suppose that $\f\in\big(B^n_{\be1}\big)_+$. Then the function {\em\rf{opf}} has $n$th derivative
in the norm
$$
\frac{d^n}{ds^n}\f(T_s)\Big|_{s=t}=n!
\underbrace{\int\!\cdots\!\int}_{n+1}(\dg^{n}\f)(\z_1,\cdots,\z_{n+1})\,d\E_t(\z_1)\,(R-T)\cdots(R-T)\,d\E_t(\z_{n+1}).
$$
\end{thm}

We refer the reader to \cite{Pe9} for proofs.

\

\section{\bf Operator H\"older--Zygmund functions}
\setcounter{equation}{0}
\label{s10}

\

As we have mentioned in the introduction, surprisingly, H\"older functions of order $\a$ must also be operator H\"older functions of order $\a$. The same is true for all spaces of the scale $\L_\a$ of H\"older--Zygmund classes.

Recall that the results of this section were obtained in \cite{AP1} (see also \cite{AP}).

Let us consider the case of unitary operators.

\begin{thm}
\label{uH}
Let $0<\a<1$ and $\f\in\L_\a$. If $U$ and $V$ are unitary operators on Hilbert space, then 
$$
\|\f(U)-\f(V)\|\le \const\|\f\|_{\L_\a}\cdot\|U-V\|^\a.
$$
\end{thm}

Note that the constant on the right-hand side of the inequality depends on $\a$.

In the proof of Theorem \ref{uH} we are going to use the following norm on $\L_\a$ (see \S 4):
$$
\|f\|_{\L_\a}=\sup_{n\ge0}2^{n\a}\|\f*W_n\|_{L^\be}+\sup_{n>0}2^{n\a}\|\f*W^\#_n\|_{L^\be}.
$$

\medskip

{\bf Proof of Theorem \ref{uH}.}
Let $\f\in\L_\a$. We have
$$
\f=\pp_+\f+\pp_-\f=\f_++\f_-.
$$
We estimate $\|\f_+(U)-\f_+(V)\|$. The norm of  $\f_-(U)-\f_-(V)$ can be estimated in the same way.
Thus we assume that $\f=\f_+$. Let
$$
\f_n\df \f*W_n.
$$
Then
\bay
\label{dia}
\f=\sum_{n\ge0}\f_n.
\ey
Clearly, we may assume that $U\ne V$. Let $N$ be the nonnegative integer such that
\bay
\label{N}
2^{-N}<\|U-V\|\le2^{-N+1}.
\ey
We have 
$$
\f(U)-\f(V)=\sum_{n\le N}\big(\f_n(U)-\f_n(V)\big)+\sum_{n> N}\big(\f_n(U)-\f_n(V)\big).
$$
By Lemma \ref{cl},
\begin{align*}
\left\|\sum_{n\le N}\big(\f_n(U)-\f_n(V)\big)\right\|&\le\sum_{n\le N}\big\|\f_n(U)-\f_n(V)\big\|\\[.2cm]
&\le\const\sum_{n\le N}2^n\|U-V\|\cdot\|\f_n\|_{L^\be}\\[.2cm]
&\le\const\|U-V\|\sum_{n\le N}2^n2^{-n\a}\|\f\|_{\L_\a}\\[.2cm]
&\le\const\|U-V\|2^{N(1-\a)}\|\f\|_{\L_\a}\le\const\|U-V\|^\a\|\f\|_{\L_\a},
\end{align*}
the last inequality being a consequence of \rf{N}.

On the other hand,
\begin{align*}
\left\|\sum_{n>N}\big(\f_n(U)-\f_n(V)\big)\right\|&\le\sum_{n>N}2\|\f_n\|_{L^\be}
\le\const\sum_{n>N}2^{-n\a}\|\f\|_{\L_\a}\\[.2cm]
&\le\const2^{-N\a}\|\f\|_{\L_\a}\le\const\|U-V\|^\a\|\f\|_{\L_\a}.\quad\bl
\end{align*}

Consider now the case of an arbitrary positive $\a$.

\begin{thm}
\label{hou}
Let $n$ be a positive integer, $0<\a<n$, and let $f\in\L_\a$. Then
for a unitary operator $U$ and a bounded self-adjoint operator $A$ on Hilbert space the following inequality holds:
$$
\left\|\sum_{k=0}^n(-1)^k\left(\begin{matrix}n\\k\end{matrix}\right)\f\big(e^{{\rm i}kA}U\big)
\right\|\le\const\|\f\|_{\L_\a}\|A\|^\a.
$$
\end{thm}

The proof of Theorem \ref{hou} given in \cite{AP1} is based on multiple operator integrals and Lemma \ref{rrp}. Let me explain
how we arrive at triple operator integrals in the case $n=2$. In this case the following formula holds:
\begin{align*}
f(U_1)-2f(U_2)+f(U_3)&=\!2\!\iiint\!(\cd^2f)(\z,\t,\up)\,dE_1(\z)(U_1-U_2)\,dE_2(\t)(U_2-U_3)\,dE_3(\up)\\[.2cm]
&+\iint(\cd f)(\z,\t)\,dE_1(\z)(U_1-2U_2+U_3)\,dE_3(\t),
\end{align*}
where $U_1$, $U_2$, and $U_3$ are unitary operators and $f$ is a function on $\T$ such that 
the function $\cd^2f$ belongs to the space $C(\T)\hat\oplus_{\rm i}C(\T)\hat\oplus_{\rm i}C(\T)$.
We refer the reader to \cite{AP1} for the proofs.

Consider now more general classes of functions. Suppose that $\o$ is a modulus of continuity, i.e.,
$\o$ is a nonnegative continuous function on $[0,\be)$ such that $\o(0)=0$ and
$$
\o(x+y)\le\o(x)+\o(y),\quad x,\,y\ge0.
$$
The class $\L_\o$ consists, by definition, of functions $\f$ such that
$$
|\f(\z)-\f(\t)|\le\const\o(|\z-\t|),\quad\z,\,\t\in\T.
$$
We put
$$
\|\f\|_{\L_\o}\df\sup_{\z\ne\t}\frac{|\f(\z)-\f(\t)|}{\o(|\z-\t|)}
$$

Given a modulus of continuity $\o$, we define
$$
\o^*(x)=x\int_x^\be\frac{\o(t)}{t^2}\,dt.
$$

\begin{thm}
\label{mco}
Let $\o$ be a modulus of continuity and let $U$ and $V$ be unitary operators on Hilbert space. Then for a function 
$\f\in\L_\o$,
$$
\|\f(U)-\f(V)\|\le\const\|\f\|_{\L_\o}\o^*\big(\|U-V\|\big).
$$
\end{thm}

Note that if $\o$ is a modulus of continuity, for which
$$
\o^*(x)\le\const\o(x),
$$
then for unitary operators $U$ and $V$ the following inequality holds:
$$
\|\f(U)-\f(V)\|\le\const\|\f\|_{\L_\o}\o\big(\|U-V\|\big).
$$

We refer the reader to \cite{AP1} for an analog of Theorem \ref{mco} for higher order moduli of continuity.

Finally, to conclude this section, I would like to mention that similar results also hold for self-adjoint operators and for contractions.
In particular, the analog of Theorem \ref{mco} for self-adjoint operators improves the estimate obtained in \cite{FN}.
We refer the reader to \cite{AP1} for detailed results.

\

\section{\bf Livshits--Krein trace formulae}
\setcounter{equation}{0}
\label{s11}

\

The spectral shift function for a trace class perturbation of a self-adjoint operator
was introduced in a special case by
I.M. Lifshitz \cite{L} and in the general case by M.G. Krein \cite{Kr}.
%It plays a very important role in perturbation theory. 
It was shown in   \cite{Kr} that for a pair of self-adjoint (not necessarily bounded) operators $A$ and $B$ 
satisfying $B-A\in\bS_1$, there exists a unique function $\xi\in L^1(\R)$ such that
\bay
\label{sa}
\trace\big(\f(B)-\f(A)\big)=\int_\R\f'(x)\xi(x)\,dx,
\ey
whenever $\f$ is a function on $\R$ such that the Fourier transform of $\f'$ is in $L^1(\R)$. The function $\xi$ is
called the {\it spectral shift function corresponding to the pair} $(A,B)$. 

A similar result was obtained in \cite{Kr2} for pairs of unitary operators $(U,V)$ with 
\lb$V-U\in\bS_1$. For each such pair there exists a function $\xi$ on the unit circle $\T$
of class $L^1(\T)$ such that 
\bay
\label{u}
\trace\big(\f(V)-\f(U)\big)=\int_\T\f'(\z)\xi(\z)\,d\m(\z),
\ey
whenever $\f'$ has absolutely convergent Fourier series. Such a function $\xi$ is unique modulo an additive constant
and it is called a {\it spectral shift function corresponding to the pair} $(U,V)$. We refer the reader
to the lectures of M.G. Krein \cite{Kr3}, in which the above results were discussed in detail
(see also \cite{BS} and the survey article \cite{BY}).

Note that the spectral shift function plays an important role in perturbation theory. We mention here
the paper \cite{BK}, in which the following remarkable formula was found:
$$
\det S(x)=e^{-2\pi{\rm i}\xi(x)},
$$
where $S$ is the scattering matrix corresponding to the pair $(A,B)$. 

It was shown later
in \cite{BS3} that formulae \rf{sa} and \rf{u} hold under less restrictive assumptions on $\f$.

Note that the right-hand sides of \rf{sa} and \rf{u} make sense for an arbitrary Lipschitz function $\f$.
However, it turns out that the condition $\f\in\Li$ (i.e., $\f$ is a Lipschitz function) does not imply
that $\f(B)-\f(A)$ or $\f(V)-\f(U)$ belongs to $\bS_1$. This is not even true for bounded $A$ and $B$ and
continuously differentiable $\f$. The first such examples were given in \cite{F1}. 

In this section we present results of \cite{Pe1} and \cite{Pe4} that give necessary conditions and sufficient conditions
on the function $\f$ for trace formulae \rf{sa} and \rf{u} to hold.

We start with the case of unitary operators. Recall that the class $\cL$ of functions on $\T$ was defined in \S\,7.

\begin{thm}
\label{nuu}
Let \mbox{$\f\in C^1(\T)$}. Suppose that \mbox{$\f\not\in\cL$}. Then there exist unitary operators $U$ and $V$ such that
$$
U-V\in\bS_1,
$$
but
$$
\f(U)-\f(V)\not\in\bS_1.
$$
\end{thm}

\begin{cor}
\label{slu}
Let \mbox{$\f\in C^1(\T)\setminus B_1^1$}. Then there exist unitary operators $U$ and $V$ such that
$$
U-V\in\bS_1,
$$
but
$$
\f(U)-\f(V)\not\in\bS_1.
$$
\end{cor}

{\bf Proof of Theorem \ref{nuu}.} As we have discussed in \S\,7, if $\f\not\in\cL$, then the divided difference $\dg\f$ is not a multiplier of the class $\frak N$ of kernel functions of trace class integral operators on $L^2(\m)$. Now, the same reasoning as in the proof of Theorem \ref{ncL} allows us to construct sequences of unitary operators $\{U_n\}_{n\ge1}$ and $\{V_n\}_{n\ge1}$ such that
$$
\lim_{n\to\be}\|U_n-V_n\|_{\bS_1}=0
$$
but
$$
\lim_{n\to\be}\frac{\|\f(U_n)-\f(V_n)\|_{\bS_1}}{\|U_n-V_n\|_{\bS_1}}=\be.
$$
It is easy to see now that we can select certain terms of these sequences with repetition (if necessary) and obtain sequences 
 $\{{\mathcal U}_n\}_{n\ge1}$ and $\{\V_n\}_{n\ge1}$
of unitary operators such that
$$
\sum_{n\ge1}\|{\mathcal U}_n-\V_n\|_{\bS_1}<\be, 
$$
but
$$
\sum_{n\ge1}\|\f({\mathcal U}_n)-\f(\V_n)\|_{\bS_1}=\be.
$$
Now it remains to define unitary operators $U$ and $V$ by
$$
U=\sum_{n\ge1}\oplus \,{\mathcal U}_n\quad\mbox{and}\quad
V=\sum_{n\ge1}\oplus \V_n.\quad\bl
$$

The following sufficient condition improves earlier results in \cite{BS3}. 

\begin{thm}
\label{duu}
Suppose that $\f\in B_{\be1}^1$. Let $U$ and $V$ be unitary operators such that $V-U\in\bS_1$ and let
$\xi$ be a spectral shift  function corresponding to the pair $(U,V)$. Then
\bay
\label{vS1}
\f(V)-\f(U)\in\bS_1
\ey
and
trace formula {\em\rf{u}} holds.
\end{thm}

\Pf Let us first prove \rf{vS1}. By Theorem \ref{du}, \mbox{$\dg\f\in C(\T)\hat\otimes C(\T)$} which implies that
$\dg\f\in {\frak M}(E_V,E_U)$. Thus by the Birman--Solomyak formula,
$$
\f(V)-\f(U)=\iint\limits_{\T\times\T}\frac{\f(\z)-\f(\t)}{\z-\t}\,dE_V(\z)\,(U-V)\,dE_U(\t).
$$
It follows that $\f(V)-\f(U)\in\bS_1$.

To prove that \rf{u} holds, we recall that by the results of \cite{Kr2},
\rf{u} holds for trigonometric polynomials. It suffices now to approximate
$\f$ by trigonometric polynomials in the norm of $B^1_{\be1}$.  $\bl$

Let us proceed to the case of self-adjoint operators. The following results were obtained in \cite{Pe4}.

\begin{thm}
\label{nus}
Suppose that $\f$ is a continuously differentiable function on $\R$ such that $\f\not\in\cL(\R)$.
Then there exist self-adjoint operators $A$ and $B$ such that 
$$
B-A\in\bS_1,
$$
but
$$
\f(B)-\f(A)\not\in\bS_1.
$$
\end{thm}

In particular, Theorem \ref{nus} implies that the condition that $\f\in B_1^1(\R)$ locally is a necessary condition for 
trace formula \rf{sa} to hold.

\begin{thm}
\label{dus}
Suppose that $\f\in B^1_{\be1}(\R)$. Let $A$ and $B$ be self-adjoint operators (not necessarily bounded) such that 
$B-A\in\bS_1$ and let $\xi$ be the spectral shift function that corresponds to the pair $(A,B)$. Then
$\f(B)-\f(A)\in\bS_1$ and trace formula {\em\rf{sa}} holds.
\end{thm}

The proof of Theorem \ref{dus} is more complicated than the proof of Theorem \ref{duu}, because nice functions are not dense
in $B_{\be1}^1(\R)$, and to prove \rf{sa} we have to use a weak approximation, see \cite{Pe4}.

\

\section{\bf Koplienko--Neidhardt trace formulae}
\setcounter{equation}{0}
\label{s12}

\

In this section we consider trace formulae in the case of perturbations of Hilbert--Schmidt class $\bS_2$.

 Let $A$ and $B$ be  self-adjoint operators such that $K\df B-A\in\bS_2$. In this case the
operator $\f(B)-\f(A)$ does not have to be in $\bS_1$ even for very nice functions $\f$. The idea of
Koplienko in \cite{Ko} was to consider the operator
$$
\f(B)-\f(A)-\frac{d}{ds}\Big(\f(A+sK)\Big)\Big|_{s=0}
$$
and find a trace formula under certain assumptions on $\f$. It was shown in \cite{Ko} that there exists
a unique function $\eta\in L^1(\R)$ such that
\bay
\label{ktf}
\trace\left(\f(B)-\f(A)-\frac{d}{ds}\Big(\f(A+sK)\Big)\Big|_{s=0}\right)=\int_\R\f''(x)\eta(x)\,dx
\ey
for rational functions $\f$ with poles off $\R$. The function $\eta$ is called the
{\it generalized spectral shift function corresponding to the pair} $(A,B)$. 

A similar problem for unitary operators was considered by Neidhardt in \cite{Ne}. Let $U$ and $V$ be
unitary operators such that $V-U\in\bS_2$. Then $V=\exp({\rm i}A)U$, where $A$ is a self-adjoint operator
in $\bS_2$. Put $U_s=e^{{\rm i}sA}U$, $s\in\R$. It was shown in \cite{Ne} that there exists a function
$\eta\in L^1(\T)$ such that
\bay
\label{ntf}
\trace\left(\f(V)-\f(U)-\frac{d}{ds}\Big(\f(U_s)\Big)\Big|_{s=0}\right)=
\int_\T\f''\eta\,d\m,
\ey
whenever $\f''$ has absolutely convergent Fourier series. Such a function $\eta$ is unique modulo a
constant and it is called a {\it generalized spectral shift function corresponding to the pair} $(U,V)$.

We state in this section results of \cite{Pe7} that guarantee the validity of trace formulae 
\rf{ktf} and \rf{ntf} under considerably less restrictive assumptions on $\f$.

\begin{thm}
\label{uk}
Suppose that $U$ and $V=e^{{\rm i}A}U$ are unitary operators on Hilbert space such that $U-V\in\bS_2$.
Let $\f\in B^2_{\be1}$. Then
$$
\f(V)-\f(U)-\frac{d}{ds}\Big(\f\big(e^{{\rm i}sA}U\big)\Big)\Big|_{s=0}\in\bS_1
$$
and trace formula {\em\rf{ntf}} holds.
\end{thm}

\begin{thm}
\label{sk}
Suppose that $A$ and $B$ are self-adjoint operators (not necessarily bounded) on Hilbert space
such that $K=B-A\in\bS_2$.
Let $\f\in B^2_{\be1}(\R)$. Then
$$
\f(B)-\f(A)-\frac{d}{ds}\Big(\f(A+sK)\Big)\in\bS_1
$$
and trace formula {\em\rf{ktf}} holds.
\end{thm}

\

\section{\bf Perturbations of class $\bS_p$}
\setcounter{equation}{0}
\label{s13}

\

In the final section of this survey article we consider the problem of the behavior of the function of an operator under perturbations by operators of 
Schatten--von Neumann class $\bS_p$. In \S\,11 we have already considered the special case of perturbations of trace class. We have seen that the condition $\f\in{\rm Lip}$ (i.e., $\f$ is a Lipschitz function) does not guarantee 
that trace class perturbations of an operator lead to trace class changes of the function of the operator. 

On the other hand, Theorem \ref{BST2} shows that for a Lipschitz function $\f$ the condition $A-B\in\bS_2$ implies that
$\f(A)-\f(B)\in\bS_2$. 

The same problem for Lipschitz functions and perturbations of class $\bS_p$, $1<p<\be$, $p\ne2$, is still open.

In the case $p<1$ the following results were found in \cite{Pe2}:

\begin{thm}
\label{dp}
Let $0<p<1$ and let $\f\in B_{\be p}^{1/p}$. Suppose that $U$ and $V$ are unitary operators such that 
$U-V\in\bS_p$. Then $\f(U)-\f(V)\in\bS_p$.
\end{thm}

\begin{thm}
\label{np}
Let $0<p<1$. Suppose that $\f$ is a continuously differentiable function on $\T$ such that
$\f(U)-\f(V)\in\bS_p$, whenever $U$ and $V$ are unitary operators such that $U-V\in\bS_p$.
Then $\f\in B_p^{1/p}$.
\end{thm}

As in the case $p=1$, Theorem \ref{np} can be improved: under the hypotheses of Theorem \ref{np},
the Hankel operators $H_\f$ and $H_{\bar\f}$ must map the Hardy class $H^1$ into the Besov space $B_p^{1/p}$.

The same results also hold for contractions and analogs of these results can also be obtained for self-adjoint operators
(in the analog of Theorem \ref{np} for self-adjoint operators the conclusion is that $\f$ belongs to $B_p^{1/p}$ locally).

We proceed now to the results of \cite{AP2} (see also \cite{AP}) that describe the behavior of $\f(U)$ for functions of class $\L_\a$ under perturbations 
of $U$ by operators of class $\bS_p$.

\medskip

{\bf Definition.} Let $p>0$. We say that a compact operator $T$ belongs to the ideal $\bS_{p,\be}$ if its singular values 
$s_n(T)$ satisfies the estimate:
$$
\|T\|_{\bS_{p,\be}}\df\sup_{n\ge0}s_n(T)(1+n)^{1/p}<\be.
$$
Clearly,
$$
\bS_p\subset\bS_{p,\be}\subset\bS_q
$$
for any $q>p$. Note that $\|\cdot\|_{\bS_{p,\be}}$ is not a norm, though for $p>1$, the space $\bS_{p,\be}$ has a norm equivalent to
$\|\cdot\|_{\bS_{p,\be}}$.

\begin{thm}
\label{apw}
Let $p\ge1$, $0<\a<1$, and let $\f\in\L_\a$. Suppose that $U$ and $V$ are unitary operators on Hilbert space such that $U-V\in\bS_p$.
Then
$$
\f(U)-\f(V)\in\bS_{\frac p\a,\be}
$$
and
$$
\|\f(U)-\f(V)\|_{\bS_{\frac p\a,\be}}\le\const\|f\|_{\L_\a}\|B-A\|_{\bS_p}^\a.
$$
\end{thm}

In the case when $p>1$ Theorem \ref{apw} can be improved by using interpolation arguments.

\begin{thm}
\label{ap}
Let $p>1$, $0<\a<1$, and let $\f\in\L_\a$. Suppose that $U$ and $V$ are unitary operators on Hilbert space such that $U-V\in\bS_p$.
Then
$$
\f(U)-\f(V)\in\bS_p
$$
and
$$
\|\f(U)-\f(V)\|_{\bS_p}\le\const\|f\|_{\L_\a}\|B-A\|_{\bS_p}^\a.
$$
\end{thm}

Note that the constants in the above inequalities depend on $\a$.

Let us sketch the proof of Theorem \ref{apw}. We refer the reader to \cite{AP2} for a detailed proof.

As in the proof of Theorem \ref{uH}, we assume that $\f\in\big(\L_\a\big)_+$ and
we consider the expansion \rf{dia}. Put
$$
Q_N=\sum_{n\le N}\big(\f_n(U)-\f_n(V)\big)\quad\mbox{and}\quad R_N=\sum_{n>N}\big(\f_n(U)-\f_n(V)\big).
$$
Then
$$
\|R_N\|\le2\sum_{n\ge N}\|\f_n\|_{L^\be}\le\const2^{-\a N}\|\f\|_{\L_\a}.
$$
It follows from Lemma \ref{cl} that 
$$
\|\f_n(U)-\f_n(V)\|_{\bS_p}\le\const2^n\|\f_n\|_{L^\be}\|U-V\|_{\bS_p}
$$
which implies that
$$
\|Q_N\|_{\bS_p}\le\const2^{(1-\a)N}\|\f\|_{\L_\a}\|U-V\|_{\bS_p}.
$$
The proof can easily be completed on the basis of the following estimates:
$$
s_n(Q_N)\le(1+n)^{-1/p}\|Q_N\|_{\bS_p}
$$
and
$$
s_n\big(\f(U)-\f(V)\big)\le s_n(Q_N)+\|R_N\|.\quad\bl
$$

Consider now the case of higher order differences.

\begin{thm}
\label{hopw}
Let $0<\a<n$ and $p\ge n$. Suppose that $U$ is a unitary operator and $A$ is a self-adjoint operator of class $\bS_p$.
Then
$$
\sum_{k=0}^n(-1)^k\left(\begin{matrix}n\\k\end{matrix}\right)\f\big(e^{{\rm i}kA}U\big)\in\bS_{\frac p\a,\be}
$$
and
$$
\left\|\sum_{k=0}^n(-1)^k\left(\begin{matrix}n\\k\end{matrix}\right)\f\big(e^{{\rm i}kA}U\big)\right\|_{\bS_{\frac p\a,\be}}
\le\const\|f\|_{\L_\a}\|A\|_{\bS_p}^\a.
$$
\end{thm}

Again, if $p>n$, Theorem \ref{hopw} can be improved by using interpolation arguments.

\begin{thm}
\label{hop}
Let $0<\a<n$ and $p>n$. Suppose that $U$ is a unitary operator and $A$ is a self-adjoint operator of class $\bS_p$.
Then
$$
\sum_{k=0}^n(-1)^k\left(\begin{matrix}n\\k\end{matrix}\right)\f\big(e^{{\rm i}kA}U\big)\in\bS_{\frac p\a}
$$
and
$$
\left\|\sum_{k=0}^n(-1)^k\left(\begin{matrix}n\\k\end{matrix}\right)\f\big(e^{{\rm i}kA}U\big)\right\|_{\bS_{\frac p\a}}
\le\const\|f\|_{\L_\a}\|A\|_{\bS_p}^\a.
$$
\end{thm}

We refer the reader to \cite{AP2} for the proofs of Theorems  \ref{hopw} and \ref{hop}.

Note that similar results also hold for contractions and for self-adjoint operators.

\

\

\noindent
\begin{tabular}{p{8cm}p{14cm}}
Department of Mathematics \\
Michigan State University  \\
East Lansing, Michigan 48824\\
USA
\end{tabular}


\begin{thebibliography}{99}
\label{bibl}

\bibitem[AP1]{AP}{\sc A.B. Aleksandrov} and {V.V. Peller}, {\em Functions  of perturbed operators}, C. R. Acad. Sci. Paris S\'er. I Math.

\bibitem[AP2]{AP1} {\sc A.B. Aleksandrov} and {V.V. Peller}, {\em Operator H\"older--Zygmund functions}, to appear.

\bibitem[AP3]{AP2} {\sc A.B. Aleksandrov} and {V.V. Peller},  {\em The behavior of functions of operators under perturbations of
class $\bS_p$}, to appear.


\bibitem[ABF]{ABF} {\sc J. Arazy, T. Barton}, and {\sc Y. Friedman}, 
{\em Operator differentiable functions},
Int. Equat. Oper. Theory  {\bf 13} (1990), 462--487.

\bibitem[ACDS] {ACDS} {\sc N.A. Azamov, A.L. Carey, P.G. Dodds,} and {\sc F.A. Sukochev},
{\it Operator integrals, spectral shift and spectral flow}, arXiv:math/0703442.

\bibitem[Be]{Be} {\sc G. Bennett}, {\em Schur multipliers}, Duke Math. J. {\bf44} (1977), 603--639.


\bibitem[BK]{BK} {\sc M.S. Birman} and {\sc M.G. Krein}, {\em On the theory of wave operators and
scattering operators}, Dokl. Akad. Nauk SSSR {\bf144} (1962), 475--478.

English transl.: Sov. Math. Dokl. {\bf3} (1962), 740--744.

\bibitem[BS1]{BS1} {\sc M.S. Birman} and {\sc M.Z. Solomyak}, {\em Double Stieltjes operator
integrals}, Problems of Math. Phys., Leningrad. Univ. {\bf1} (1966), 33--67 (Russian).

English transl.: Topics Math. Physics {\bf1} (1967), 25--54, Consultants Bureau Plenum
Publishing Corporation, New York.

\bibitem[BS2]{BS2} {\sc M.S. Birman} and {\sc M.Z. Solomyak}, {\em Double Stieltjes operator
integrals. II}, Problems of Math. Phys., Leningrad. Univ. {\bf2} (1967), 26--60 (Russian).

English transl.: Topics Math. Physics {\bf2} (1968), 19--46, Consultants Bureau Plenum
Publishing Corporation, New York.

\bibitem[BS3]{BS} {\sc M.S. Birman} and {\sc M.Z. Solomyak}, {\em Remarks on the spectral
shift function},  Zapiski Nauchn. Semin. LOMI {\bf27} (1972), 33--46 (Russian).

English transl.: J. Soviet Math. {\bf3} (1975), 408--419.

\bibitem[BS4]{BS3} {\sc M.S. Birman} and {\sc M.Z. Solomyak}, {\em Double Stieltjes operator
integrals. III}, Problems of Math. Phys., Leningrad. Univ. {\bf6} (1973), 27--53 (Russian).

\bibitem[BS5]{BS4} {\sc M.S. Birman} and {\sc M.Z. Solomyak},
{\it Tensor product of a finite number of spectral measures is always a spectral measure},
Integral Equations Operator Theory {\bf24} (1996), 179--187.

\bibitem[BS6]{BS5} {\sc M.S. Birman} and {\sc M.Z. Solomyak}, {\em Double operator integrals
in Hilbert space},  Int. Equat. Oper. Theory  {\bf47}  (2003), 131--168.

\bibitem[BY]{BY} {\sc M.S. Birman} and {\sc D.R. Yafaev}, {\em The spectral shift function. The
papers of M. G. Kre\u\i n and their further development},  Algebra i Analiz  {\bf4} 
(1992), 1--44 (Russian).  

English transl.: St. Petersburg Math. J.  {\bf4}  (1993), 833--870.


\bibitem[DK]{DK} {\sc Yu.L. Daletskii} and {S.G. Krein}, {\it Integration and differentiation of
functions of Hermitian operators and application to the theory of perturbations} (Russian), Trudy Sem.
Functsion. Anal., Voronezh. Gos. Univ. {\bf1} (1956), 81--105.

\bibitem[Fa1]{F0}  {\sc Yu.B. Farforovskaya}, {\em  The connection of the Kantorovich-Rubinshtein metric for spectral resolutions of selfadjoint operators with functions of operators},
Vestnik Leningrad. Univ.  {\bf19}  (1968), 94--97. (Russian).

\bibitem[Fa2]{F1}  {\sc Yu.B. Farforovskaya}, {\em An example of a Lipschitzian function of selfadjoint
operators that yields a nonnuclear increase under a nuclear perturbation}.  Zap. Nauchn. Sem.
Leningrad. Otdel. Mat. Inst. Steklov. (LOMI)  {\bf30}  (1972), 146--153 (Russian).

\bibitem[Fa3]{F}  {\sc Yu.B. Farforovskaya}, {\em An estimate of the norm of $\mid f(B)-f(A)\mid $ for
selfadjoint operators $A$ and $B$},   Zap. Nauchn. Sem. Leningrad. Otdel. Mat. Inst.
Steklov. (LOMI)  {\bf56}  (1976), 143--162 (Russian).

\bibitem[FN]{FN}  {\sc Yu.B. Farforovskaya} and {\sc L. Nikolskaya}, {\em Modulus of continuity of operator functions}, 
 Algebra i Analiz  {\bf20:3}  (2008),  224--242.

\bibitem[FrJ]{FrJ} {\sc M. Frazier} and {\sc B. Jawerth}, {\em A discrete transform and decompositions of distribution spaces},
 J. Funct. Anal. {\bf93}  (1990), 34--170. 

\bibitem[JTT]{JTT} {\sc K. Jushchenko, I.G. Todorov}, and {\sc L. Turowska},
{\em Multidimensional operator multipliers}, 	arXiv:math/0701645v2.

%\bibitem[dPS]{dPS} {\sc B. de Pagter} and {\sc F.A. Sukochev}, {\em Differentiation of operator
%functions in noncommutative \lb$L_p$-spaces}, J. Funct. Anal. {\bf212} (2004), 28--75.

\bibitem[KA]{Ka}  {\sc T. Kato}, {\em Continuity of the map $S\mapsto \mid S\mid $ for linear operators,}  Proc. Japan Acad.  {\bf49}  (1973), 157--160. 

\bibitem[KS] {KS} {\sc E. Kissin} and {V.S. Shulman}, {\it On fully operator Lipschitz functions},  J. Funct. Anal. 
{\bf 253}  (2007),   711--728.


\bibitem[Ko]{Ko} {\sc L.S. Koplienko}, {\em The trace formula for perturbations of nonnuclear type},
Sibirsk. Mat. Zh. {\bf25:5}  (1984), 62--71 (Russian).

English transl.: Sib. Math. J. {\bf25} (1984), 735--743.

%\bibitem[Ko2]{Ko2} {\sc L.S. Koplienko}, {\em Regularized spectral shift function for one-dimensional
%Schr\"{o}dinger operator with slowly decreasing potential},   Sibirsk. Mat. Zh.  {\bf26:3}  (1985),  
%72--77 (Russian).
%
%English transl.: Sib. Math. J. {\bf26} (1985), 365--369.

%\bibitem[LSY]{LSY} {\sc H. Langer, H. S. V. de Snoo} and {\sc V. A. Yavrian},
%{\em A relation for the spectral shift function of two self-adjoint extensions}, In: Recent advances in
%operator theory and related topics (Szeged, 1999), 437--445, Oper. Theory Adv. Appl., {\bf127},
%Birkh\"auser, Basel, 2001. 



\bibitem[Kr1]{Kr} {\sc M.G. Krein}, {\em On a trace formula in perturbation theory},
Mat. Sbornik {\bf33} (1953), 597--626 (Russian).

\bibitem[Kr2]{Kr2} {\sc M.G. Krein}, {\em On perturbation determinants and a trace formula for unitary and
self-adjoint operators},  Dokl. Akad. Nauk SSSR {\bf144} (1962) 268--271 (Russian).

English transl. in: Topics in integral and differential equations and operator theory, Birkh\"auser,
Basel, 1983, 107--172.

\bibitem[Kr3]{Kr3} {\sc M.G. Krein}, {\em On some new investigations in the perturbation theory of
self-adjoint operators}, in: The First Summer Math. School, Kiev, 1964, 103--187 (Russian).

\bibitem[L]{L} {\sc I.M. Lifshitz}, {\em On a problem in perturbation theory 
connected with quantum statistics}, Uspekhi Mat. Nauk {\bf7} (1952), 171--180 (Russian).


\bibitem[MM]{MM} {\sc M.M. Malamud} and {\sc S.M. Malamud}, {\em 
Spectral theory of operator measures in a Hilbert space},  Algebra i Analiz  {\bf15:3}  (2003), 1--77.

English transl.:  St. Petersburg Math. J.  {\bf15:3}  (2004), 323--373.

\bibitem[Na]{Na} {\sc M.A. Naimark}, {\it Spectral functions of symmetric operator},
Izvestia Akad. Nauk SSSR, Ser. Matem. {\bf4:3} (1940), 277--318 (Russian).

\bibitem[Ne]{Ne} {\sc H. Neidhardt}, {\em Spectral shift function and Hilbert--Schmidt
perturbation: extensions of some work of L.S. Koplienko}, Math. Nachr. {\bf 138} (1988), 7--25.

\bibitem[Pa]{Pa} {\sc B.S. Pavlov}, {\em  On multiple operator integrals}, 
Problems of Math. Anal., No. 2: Linear Operators and Operator Equations (Russian),  99--122.
Izdat. Leningrad. Univ., Leningrad, 1969. 

%\bibitem[Pa]{Pa} {\sc B.S. Pavlov}, {\em  On multiple operator integrals}, 
%Problems of Math. Anal., No. 2: Linear Operators and Operator Equations (Russian),  99--122.
%Izdat. Leningrad. Univ., Leningrad, 1969. 

\bibitem[Pee]{Pee} {\sc J. Peetre}, 
{\em New thoughts on Besov spaces}, Duke Univ. Press., Durham, NC, 1976.

%\bibitem[Pe1]{Pe}
%{\sc V.V.Peller}, {\em Estimates of functions of power bounded operators
%on Hilbert space}, J. Oper. Theory, {\bf 7} (1982), 341-372.

\bibitem[Pe1]{Pe0} {\sc V.V. Peller},
 {\em Hankel operators of class ${\bf S}_{p}$ 
and their applications (rational approximation, Gaussian processes, 
the problem of majorizing operators)}, Mat. Sbornik, 
{\bf 113} (1980), 538-581. 

English Transl. in Math. USSR Sbornik, {\bf 41}
(1982), 443-479.


\bibitem[Pe2]{Pe1} {\sc V.V. Peller}, {\em Hankel operators in the theory of perturbations
of unitary and self-adjoint operators},  Funktsional. Anal. i Prilozhen. {\bf19:2}  (1985), 
37--51 (Russian).

English transl.: Funct. Anal. Appl. {\bf19} (1985) , 111--123.

\bibitem[Pe3]{Pe2} {\sc V.V. Peller},
{\em For which $f$ does $A-B\in{\bf S}_{p}$ 
imply that $f(A)-f(B)\in{\bf S}_{p}$?}, Operator Theory, Birkh\"{a}user,
{\bf 24} (1987), 289-294.

%\bibitem[Pe4]{Pe3} {\sc V.V. Peller}, {\em When is a function of a Toeplitz operator close to a Toeplitz
%operator?}  Toeplitz operators and spectral function theory,  59--85, Oper. Theory Adv. Appl., 42,
%Birkh\"{a}user, Basel, 1989.

\bibitem[Pe4]{Pe4} {\sc V.V. Peller} {\em Hankel operators in the perturbation theory of 
of unbounded self-adjoint operators}.  Analysis and partial differential equations,  529--544,
Lecture Notes in Pure and Appl. Math., {\bf122}, Dekker, New York, 1990.

\bibitem[Pe5]{Pe5} {\sc V.V. Peller},
{\em Functional calculus for a pair of almost
commuting selfadjoint operators}, J. Funct. Anal., {\bf 112} (1993), 325--345. 

\bibitem[Pe6]{Pe6} {\sc V.V. Peller}, {\em Hankel operators and their applications,}
Springer-Verlag, New York, 2003.

\bibitem[Pe7]{Pe7} {\sc V.V. Peller}, {\em An extension of the Koplienko--Neidhardt trace formulae},
J. Funct. Anal. {\bf221} (2005), 456--481.

\bibitem[Pe8] {Pe8}  {\sc V.V. Peller}, {\em Multiple operator integrals and higher operator
derivatives}, J. Funct. Anal.  {\bf233}  (2006),  515--544.

\bibitem[Pe9] {Pe9}  {\sc V.V. Peller}, {\em Differentiability of functions of contractions}, to appear.

\bibitem[St]{S} {\sc V.V. Sten'kin}, {\em Multiple operator integrals}, 
Izv. Vyssh. Uchebn. Zaved. Matematika {\bf4 (79)} (1977), 102--115 (Russian).

English transl.: Soviet Math. (Iz. VUZ) {\bf21:4} (1977), 88--99.

\bibitem[SNF]{SNF} {\sc B. Sz.-Nagy and C. Foia\c s}, 
{\it Harmonic analysis of operators on Hilbert
space,} Akad\'{e}miai Kiad\'{o}, Budapest, 1970.


%\bibitem[SS]{SS} {\sc M.Z. Solomyak} and {\sc V.V. Sten'kin},
%{\em A certain class of multiple operator Stieltjes integrals},
%Problems of Math. Anal., No. 2: Linear Operators and Operator Equations (Russian), 122--134.
%Izdat. Leningrad. Univ., Leningrad, 1969. 

\bibitem[T]{T} {\sc H. Triebel}, {\em Theory of function spaces,}
Monographs in Mathematics, {\bf78}, Birkh\"auser Verlag, Basel, 1983.

\bibitem[W]{W} {\sc H. Widom}, {\it When are differentiable functions differentiable?}, In: Linear and
Complex Analysis Problem Book, Lect. Notes Math. {\bf 1043} (1984), 184--188.


\end{thebibliography}
\end{document}